\pdfoutput=1
\documentclass[11pt]{article}

\usepackage[a4paper,margin=1in]{geometry}
\usepackage{graphicx}
\usepackage{multirow}
\usepackage{array}
\usepackage{float}
\usepackage{amsmath,amssymb,amsfonts}
\usepackage{amsthm}
\usepackage[title]{appendix}
\usepackage{booktabs}
\usepackage{enumitem}
\usepackage{xcolor}
\usepackage{algorithm}
\usepackage{algorithmicx}
\usepackage{algpseudocode}
\usepackage[numbers,sort&compress]{natbib}
\usepackage[colorlinks=true,citecolor=blue,linkcolor=blue,urlcolor=blue]{hyperref}
\usepackage{cleveref}

\theoremstyle{plain}
\newtheorem{theorem}{Theorem}
\newtheorem{proposition}[theorem]{Proposition}
\newtheorem{corollary}[theorem]{Corollary}
\newtheorem{lemma}[theorem]{Lemma}

\theoremstyle{definition}
\newtheorem{assumption}{Assumption}
\newtheorem{definition}{Definition}

\theoremstyle{remark}

\def\x{\mathbf{x}}
\def\y{\mathbf{y}}
\def\cS{\mathcal{S}}
\def\cM{\mathcal{M}}
\def\cN{\mathcal{N}}

\def\Hess{\mathcal{H}}
\def\RBDA{\ensuremath{\mathrm{RBDA}}}
\def\RBDAnoBB{\RBDA{} w/o BB}

\def\Rm{\mathbb{R}^m}

\def\Sopt{\cS^{\mathrm{opt}}}
\def\Exp{\operatorname{Exp}}
\def\Rn{\mathbb{R}^n}
\def\cG{\mathcal{G}}
\def\cU{\mathcal{U}}

\def\cT{\mathcal{T}}
\def\cR{\mathcal{R}}
\def\br{\mathbb{R}}
\def\txm{\cT_\x\cM}
\def\rd{\mathrm{d}}
\def\D{\mathrm{D}}
\def\vp{\varphi}

\def\ak{\alpha_k}
\def\cX{\mathcal{X}}

\def\by{\bar{\mathbf{y}}}
\def\bB{\bar{B}}

\def\ty{\tilde{\mathbf{y}}}

\def\St{\mathrm{St}}
\def\tr{\text{tr}}

\begin{document}

\title{Riemannian Bilevel Optimization with Gradient Aggregation}

\author{
  Zhuo Chen\textsuperscript{1},
  Xinjian Xu\textsuperscript{2},
  Shihui Ying\textsuperscript{3}\thanks{Corresponding author:
  \href{mailto:shying@shu.edu.cn}{shying@shu.edu.cn}},
  and Tieyong Zeng\textsuperscript{4}\\[0.75em]
  \small\textsuperscript{1}Department of Mathematics, Shanghai University, Shanghai, China\\
  \small\textsuperscript{2}Qian Weichang College, Shanghai University, Shanghai, China\\
  \small\textsuperscript{3}School of Mechanics and Engineering Science, Shanghai University, Shanghai, China\\
  \small\textsuperscript{4}Department of Mathematics, The Chinese University of Hong Kong, Hong Kong, China\\[0.5em]
  \small\textsuperscript{1}\href{mailto:98chenzhuo@shu.edu.cn}{98chenzhuo@shu.edu.cn};
  \textsuperscript{2}\href{mailto:xinjxu@shu.edu.cn}{xinjxu@shu.edu.cn};
  \textsuperscript{3}\href{mailto:shying@shu.edu.cn}{shying@shu.edu.cn};
  \textsuperscript{4}\href{mailto:zeng@math.cuhk.edu.hk}{zeng@math.cuhk.edu.hk}
}
\date{}

\maketitle

\begin{abstract}
We study bilevel optimization on Riemannian manifolds when the lower-level solution set is a positive-dimensional submanifold, so that implicit differentiation fails.
We propose Riemannian Bilevel Descent Aggregation (RBDA), which extends bilevel descent aggregation to manifolds.
Its inner loop aggregates the lower-level descent direction with the upper-level gradient under a decaying multiplier, and its hypergradient is the reverse-mode derivative of the unrolled loop.
Under geodesic convexity and quadratic growth of the lower level, the inner iterates converge to a point of the optimistic solution set at a polynomial rate.
Approximate minimizers of the objective with a finite number of inner iterations converge to minimizers of the optimistic value.
In the experiments RBDA selects the optimistic solution where the implicit and unrolled estimators remain at the initial point or stop at a larger query loss.
While each of its outer steps costs more than that of the unrolled estimator, it attains the highest test accuracy in data hyper-cleaning.
\end{abstract}

{\raggedright\noindent\textbf{Keywords:} \mbox{Riemannian Optimization}, \mbox{Bilevel Optimization}, \mbox{Gradient Aggregation}, \mbox{Optimistic Solution}, \mbox{Hadamard Manifolds}\par}


\section{Introduction}\label{sec:intro}
\emph{Bilevel optimization} (BLO) provides a principled framework for tackling hierarchical decision-making problems that arise in many domains characterized by nested structures~\cite{beck2023survey,lu2024first-order}.
It has proven to be a powerful tool in a wide range of machine learning applications, including hyperparameter optimization~\cite{franceschi2018bilevel,shaban2019truncated}, meta-learning~\cite{rajeswaran2019meta,ji2020convergence}, reinforcement learning~\cite{gammelli2023graph,liu2021investigating}, and generative adversarial networks~\cite{goodfellow2014generative,metz2017unrolled}.
A general (optimistic) bilevel optimization\footnote{This work focuses on the optimistic formulation of bilevel optimization.
For pessimistic bilevel optimization, which is closely related to min--max problems, Pareto optimality, and Nash equilibrium, we refer the interested reader to~\cite{gibbons1992primer,beck2023survey}.} takes the following formulation:
$$
\min_{\x\in\Rm, \y\in\Rn} F(\x,\y)
\quad\text{s.t.}\quad \y\in\cS(\x):=\arg\min_\y f(\x,\y),
$$
where the \emph{upper-level} (UL) and \emph{lower-level} (LL) objective functions $F$ and $f$ are both jointly continuously differentiable with respect to the UL variable $\x$ and the LL variable $\y$.
The set-valued mapping $\cS(\cdot)$ indicates the solution set of the LL subproblem parameterized by $\x$.

In many applications the variables are subject to structural constraints that Euclidean space does not capture.
Examples include orthogonality~\cite{edelman1998geometry}, low-rankness~\cite{vandereycken2013lowrank}, and invariance under certain transformations~\cite{absil2009optimization}.
Many of these constraint sets are nonconvex, and Euclidean methods that accept them do not exploit their geometry~\cite{dutta2024riemannian}.
Each of them is a Riemannian manifold, and the manifold formulation keeps every iterate feasible.
It can also restore convexity, since a function that is nonconvex in the Euclidean sense may be convex along the geodesics of the manifold~\cite{sra2015conic}.
The Karcher mean of symmetric positive definite (SPD) matrices is the standard example~\cite{zhang2016first-order}.
This is the setting of \emph{Riemannian bilevel optimization} (RBLO), in which the variables of one or both levels lie on manifolds.

We now present the formal definition of the RBLO problem, in which both the UL and LL variables are constrained to lie on Riemannian manifolds:
\begin{equation}\label{P:RBLO1}
\min_{\x\in\cM,\y\in\cN} F(\x,\y)\quad\text{s.t.}\quad\y\in\cS(\x)=\arg\min_{\y\in\cN}f(\x,\y),
\end{equation}
where $\cM$ and $\cN$ are connected, complete Riemannian manifolds,%
\footnote{Euclidean space with its usual metric is itself a connected, complete Riemannian manifold, so the formulation covers the case where only one level is curved.
The few-shot and hyper-cleaning experiments run a Euclidean lower level.}
and both $F$ and $f$ are twice continuously differentiable ($C^2$) objective functions defined on the product manifold $\cM\times\cN$.

Existing RBLO methods require the \emph{Lower-Level Singleton} (LLS) condition, under which the lower level is geodesically strongly convex and $\cS(\x)$ is a single point $\y^*(\x)$.
Under this condition, the hypergradient of the value function $F(\x,\y^*(\x))$ is obtained from the Riemannian implicit function theorem~\cite{han2024framework,li2025riemannian,dutta2024riemannian}.
However, this condition fails for the lower-level problems of interest.
Symmetry orbits of the lower objective, blocks that it ignores, and rank deficiency each produce a lower solution set $\cS(\x)$ that is a positive-dimensional totally geodesic submanifold.
In this setting, the lower-level sensitivity problem is ill-posed, since the lower-level Hessian is singular in the tangent directions to $\cS(\x)$.
The standard implicit differentiation formula therefore breaks down.
Direct inversion, conjugate gradients and the Neumann series all require this inverse, so none of them helps.
Forcing uniqueness back by a strongly convex regularizer is known to be unreliable even as the regularization vanishes~\cite{liu2022general}.
This paper selects the lower solution by an inner iteration that aggregates the gradients of the two levels, and obtains the hypergradient by differentiating through that iteration.

\subsection{Related work}\label{subsec:related}
\textbf{Riemannian optimization}
Riemannian optimization treats a structured constraint set as a manifold and works intrinsically on it.
One of the earliest systematic studies in this direction was conducted by Absil et al.~\cite{absil2004riemannian}, who introduced a geometric framework for optimization algorithms on Riemannian manifolds.
In addition to that foundational work, the analysis of Riemannian optimization has been further advanced in~\cite{absil2009optimization,boumal2014optimization,sato2021riemannian,boumal2023introduction,zhang2016first-order}.
These references pay particular attention to applications such as matrix completion, machine learning, and signal processing.

\noindent
\textbf{Bilevel optimization}
Gradient-based bilevel methods in Euclidean space follow three routes.
\emph{Implicit differentiation} applies the implicit function theorem to the lower-level optimality condition and evaluates the hypergradient through Hessian--vector solves, approximated by conjugate gradients or Neumann series~\cite{pedregosa2016hyperparameter,grazzi2020iteration,ji2021bilevel}.
\emph{Iterative differentiation} unrolls a finite lower-level trajectory and backpropagates through it~\cite{franceschi2017forward,franceschi2018bilevel,shaban2019truncated}.
The iteration complexities of the two routes are compared in~\cite{grazzi2020iteration,ji2021bilevel}.
A third route reformulates the bilevel constraint through structured penalty or minimax subproblems~\cite{lu2024first-order}, in the same spirit as the Lagrangian reformulation of~\cite{dutta2024riemannian} below.
These routes say how the hypergradient is computed, and each of them is stated for a unique lower-level solution.
When the solution set is larger, a method must also choose a member of it.
The Bilevel Descent Aggregation (BDA) line of Liu, Zhang, and coauthors~\cite{liu2020generic,liu2022general} makes that choice implicitly, by aggregating the two descent directions inside the unrolled dynamic.
In this way, the inner iterates select the member of $\cS(\x)$ preferred by the upper objective, which is the \emph{optimistic} solution.
Masiha et al.~\cite{masiha2026select} make the choice explicitly, by low-temperature sampling, and then differentiate through the pseudoinverse of the singular lower-level Hessian on the normal space.
Their analysis is Euclidean, and it assumes a local Polyak--\L{}ojasiewicz condition together with a unique nondegenerate optimistic minimizer.

\noindent
\textbf{Riemannian Bilevel Optimization}
Bilevel problems with a manifold at one level include meta-learning with parameters on the Stiefel or Grassmann manifold~\cite{han2024framework,li2025riemannian}, the hyper-representation of SPD data through Stiefel embeddings~\cite{han2024framework,shi2025adaptive}, and robust estimation on SPD matrices, whose lower level is the Karcher mean~\cite{li2025riemannian,shi2025adaptive}.
Li and Ma~\cite{li2025riemannian} work out the Riemannian cross-derivative that the hypergradient requires and build the deterministic RieBO and the stochastic RieSBO algorithms on it.
Concurrently, Han et al.~\cite{han2024framework} derive the hypergradient from the Riemannian implicit function theorem and estimate it by Hessian inversion or a truncated Neumann series.
Their guarantees hold under the LLS condition, namely a geodesically strongly convex lower level with a Lipschitz Riemannian gradient.
Every estimator in the two frameworks therefore queries the lower-level Hessian.

Dutta et al.~\cite{dutta2024riemannian} work with first-order information alone.
They impose the lower-level optimality as a constraint and adjust the multiplier of the resulting Lagrangian reformulation, so that the single-loop RF$^2$SA method reaches $\epsilon$-stationarity of $\vp$.
The guarantee covers stochastic setups as well, and the scheme still presumes a well-defined lower-level solution.
Most recently, Shi et al.~\cite{shi2025adaptive} propose AdaRHD, a fully adaptive stepsize strategy that removes all prior knowledge of the problem constants and of the curvature.
It retains the $O(1/\epsilon)$ complexity of non-adaptive methods, and replacing exponential maps by retractions preserves the bound.
Their count is of outer iterations to $\epsilon$-stationarity of $\vp$ under LLS, so it is not comparable to the per-hypergradient inner counts of this paper.

Each of these methods is analyzed under the LLS condition, and none of them gives a guarantee when it fails.

\subsection{Contributions and Organization}\label{subsec:contrib}
Our contributions are as follows.
\begin{itemize}[label=--, leftmargin=*]
  \item We extend the bilevel descent aggregation of Liu et al.~\cite{liu2020generic,liu2022general} to Riemannian manifolds as Riemannian Bilevel Descent Aggregation (\RBDA{}), whose inner loop aggregates the lower-level descent direction with the upper-level gradient in the tangent space.
  The hypergradient is obtained by reverse-mode automatic differentiation through the unrolled inner loop.
  \item We prove that, without the LLS condition, the inner iterates of \RBDA{} converge to a point of the optimistic solution set when the lower level is g-convex on a Hadamard manifold and satisfies a quadratic growth condition, with rates for both objectives under a growth condition of the upper level.
  We also show that approximate minimizers of the objective with a finite number of inner iterations converge to minimizers of the optimistic value.
  For a Euclidean lower level, these results improve the lower-level rate of~\cite{liu2022general} and add an upper-level rate.
  \item We test \RBDA{} on two synthetic problems with and without the LLS condition and on few-shot meta-learning and data hyper-cleaning over three datasets.
  Under the LLS condition it performs comparably to the implicit and unrolled estimators.
  In the non-LLS case it converges to the optimistic solution, whereas the implicit and unrolled estimators do not.
  In few-shot meta-learning its meta-test accuracy exceeds that of the implicit and unrolled estimators on two of the three datasets and matches it on the third.
  In data hyper-cleaning it attains the highest test accuracy on all three datasets.
\end{itemize}

\noindent
The remainder of the paper is organized as follows.
Section~\ref{sec:preliminaries} collects the Riemannian background and the notation.
Section~\ref{sec:algorithm} proposes a Riemannian bilevel algorithm based on gradient aggregation, together with its hypergradient and stepsizes.
Section~\ref{sec:analysis} analyzes the convergence of the proposed algorithm.
Section~\ref{sec:experiments} reports the numerical experiments, and Section~\ref{sec:conclusion} concludes the paper.
Appendix~\ref{app:proofs} contains the proofs of Section~\ref{sec:analysis}, and Appendix~\ref{app:supp} contains the supplementary experiments.

\section{Preliminaries}\label{sec:preliminaries}
In this section, we provide the necessary background and notation, the latter summarized in Table~\ref{tab:symbols}.
We begin by reviewing the concepts from Riemannian optimization~\cite{boumal2023introduction,boumal2014optimization}, which serve as the geometric foundation of our framework.

Suppose $\cM$ is an $m$-dimensional differentiable manifold.
Given a differentiable curve $\gamma:[-\delta,\delta]\to\cM$ with $\gamma(0)=\x$ for some $\delta>0$, the \emph{tangent space} at $\x\in\cM$, denoted $\cT_\x\cM$, is the collection of all derivatives $\gamma^\prime(0)$ of such curves, and its elements are called \emph{tangent vectors}.
The \emph{tangent bundle} $\cT\cM:=\bigcup_{\x\in\cM}\cT_\x\cM$ collects these spaces over the entire manifold.

The notion of Riemannian manifold is defined as follows.
\begin{definition}\label{D:RieMan}
A \textbf{Riemannian manifold} is a smooth manifold $\cM$ endowed with a metric $g$.
Here, $g$ assigns a smoothly varying inner product $\langle\cdot,\cdot\rangle_\x$ to each tangent space $\cT_\x\cM$.\footnote{See~\cite{boumal2023introduction} p.~41, Definition 3.53 for more details.}
We omit the subscript of the inner product when it causes no ambiguity.
\end{definition}

The \textbf{Levi--Civita connection} $\nabla$ of $(\cM,g)$ is the unique affine connection that is torsion-free and compatible with the metric, and it provides the covariant derivative $\nabla_XY$ of a vector field $Y$ along a vector field $X$~\cite{boumal2023introduction}.
For a two-dimensional subspace $\Pi\subset\cT_\x\cM$, the \textbf{sectional curvature} $\sec(\Pi)$ is the Gaussian curvature at $\x$ of the surface swept out by the geodesics whose initial velocities lie in $\Pi$~\cite[Ch.~8]{lee2018introduction}.
Then we review the notion of the differential of a mapping between manifolds.

\begin{definition}\label{D:DifGra}
Let $F:\cM\to\cN$ be a smooth mapping between two differentiable manifolds.
At each point $\x\in\cM$, the \textbf{differential} of $F$ is the linear mapping
$$
\D F:\cT_\x\cM\to\cT_{F(\x)}\cN
$$
that sends the velocity $\gamma^\prime(0)$ of a curve $\gamma$ through $\x$ to the velocity $(F\circ\gamma)^\prime(0)$ of its image.
For a smooth function $f:\cM\to\br$, its differential is denoted by $\D f$.

The \textbf{Riemannian gradient} of $f$ at $\x\in\cM$ is defined by
$$
    \D f(\xi)=\langle\cG f,\xi\rangle, \forall\xi\in\cT_\x\cM,
$$
where $\langle\cdot,\cdot\rangle$ denotes the Riemannian metric on tangent space $\cT_\x\cM$.

The \textbf{Riemannian Hessian} of $f$ at $\x\in\cM$ is the covariant derivative of the Riemannian gradient, that is, the self-adjoint linear operator $\Hess f(\x):\cT_\x\cM\to\cT_\x\cM$ given by $\Hess f(\x)[\xi]=\nabla_\xi\,\cG f$ for every $\xi\in\cT_\x\cM$, with $\nabla$ the Levi--Civita connection.
\end{definition}

The Riemannian gradient is a natural generalization of the Euclidean gradient to Riemannian manifolds.
For a function of two manifold variables $f:\cM\times\cN\to\br$, we write $\cG_\x f$ and $\cG_\y f$ for the partial Riemannian gradients in the $\cM$- and $\cN$-directions, and we use the associated block second-order derivatives.
The \textbf{lower Hessian} $\Hess_\y f(\x,\y):\cT_\y\cN\to\cT_\y\cN$ is the covariant derivative of $\cG_\y f$ along $\cN$, and the \textbf{mixed second-order derivative} (cross-Hessian) $\Hess_{\x\y}f(\x,\y):\cT_\x\cM\to\cT_\y\cN$ is the covariant derivative of $\cG_\y f$ along $\cM$.
It measures how the lower-level gradient responds to a change in the upper variable $\x$.
Its Riemannian adjoint is denoted $\Hess_{\x\y}f(\x,\y)^{*}:\cT_\y\cN\to\cT_\x\cM$.
We next introduce the geometric tools for moving along a manifold and measuring distances on it.

\begin{definition}
A \textbf{geodesic} on a Riemannian manifold $(\cM,g)$ is a smooth curve $\gamma:I\subset\br\to\cM$, whose velocity vector field $\dot{\gamma}(t)\in \cT_{\gamma(t)}\cM$ satisfies the geodesic equation
$$
	\nabla_{\dot{\gamma}(t)}\dot{\gamma}(t)=0,
$$
where $\nabla$ denotes the Levi--Civita connection associated with metric $g$.
In Euclidean space $(\Rn, \langle\cdot,\cdot\rangle)$, this condition reduces to $\ddot{\gamma}(t)=0$, and thus a geodesic coincides with a straight line.
\end{definition}

The exponential and logarithmic maps connect a manifold with its tangent spaces along these geodesics.

\begin{definition}
For $\x\in\cM$ and $\xi\in \cT_\x\cM$, let $\gamma_\xi:[0,1]\to\cM$ be the (unique) geodesic with $\gamma_\xi(0)=\x$ and $\dot\gamma_\xi(0)=\xi$.
The \textbf{exponential map} at $\x$ is $\Exp_\x(\xi)=\gamma_\xi(1)$, the point reached by moving along that geodesic for unit time, so that $\Exp_\x$ transfers the linear structure of $\cT_\x\cM$ to the manifold.
For $\y$ in a normal neighborhood of $\x$, the \textbf{logarithmic map} is the local inverse of $\Exp_\x$, namely $\Exp^{-1}_\x(\y)\in\cT_\x\cM$ such that $\Exp_\x\bigl(\Exp^{-1}_\x(\y)\bigr)=\y$.
In a geodesically normal neighborhood, the Riemannian distance between two points $\x$ and $\y$ is given by the norm of the logarithmic map, that is, $\rd(\x,\y)=\|\Exp^{-1}_\x(\y)\|_\x$.
\end{definition}

A \textbf{Hadamard manifold} is a complete, simply connected Riemannian manifold whose sectional curvature is nonpositive everywhere.
On such a manifold the exponential map at every point is a global diffeomorphism and any two points are joined by a unique geodesic~\cite[Thm.~12.8 and Prop.~12.9]{lee2018introduction}.

Although the exponential map provides a mathematically precise means of following geodesics, it is often computationally expensive or unavailable in closed form.
Geodesics on the Stiefel manifold, for instance, involve matrix exponentials~\cite{edelman1998geometry}.
To address this, we adopt the more practical notion of a retraction, a smooth mapping from the tangent space to the manifold that locally approximates the exponential map and enables efficient updates during optimization.

\begin{definition}
A \textbf{retraction} on a manifold is a smooth mapping $\cR$ from the tangent bundle $\cT\cM$ onto $\cM$ whose restriction $\cR_\x$ to $\cT_\x\cM$ satisfies $\cR_\x(0)=\x$ and $\D\cR_\x(0)=\mathrm{Id}_{\cT_\x\cM}$ for every $\x\in\cM$.
\end{definition}

While retraction maps a tangent vector back onto the manifold, many optimization algorithms also require transporting vectors between different tangent spaces.
This is achieved through \emph{parallel transport}, which allows the comparison and combination of vectors defined at different points on the manifold in a geometrically consistent way.

\begin{definition}
Given a Riemannian manifold $\cM$ and two points $\x,\y\in\cM$, the \textbf{parallel transport} $\Gamma_\x^\y:\txm\to\cT_\y\cM$ is a linear operator that preserves the inner product, so that $\langle\Gamma_\x^\y\xi,\Gamma_\x^\y\zeta\rangle_\y=\langle\xi,\zeta\rangle_\x$ for all $\xi,\zeta\in\txm$.
\end{definition}

Next, we review the definition of geodesic (strong) convexity and geodesic Lipschitz smoothness, which are the generalizations of the Euclidean counterparts.

\begin{definition}
A \textbf{geodesically convex} neighborhood $\cU\subset\cM$ is a set such that for any two points in the set, there exists a unique geodesic connecting them that lies entirely in $\cU$.

A function $f:\cU\to\br$ is called \textbf{geodesically convex} (\textbf{g-convex})~\cite{zhang2016first-order} if for any $p, q\in\cU$, we have $f(\gamma(t))\le(1-t)f(p)+tf(q)$, where $\gamma$ is a geodesic in $\cU$ with $\gamma(0)=p$ and $\gamma(1)=q$.

If $f$ is continuously differentiable, it is g-convex if and only if
$$
f(\y)\geq f(\x)+\langle\cG f(\x),\Exp^{-1}_\x(\y)\rangle_\x ,
$$
and it is \textbf{$C_f$-strongly g-convex} if and only if
$$
f(\y)\geq f(\x)+\langle\cG f(\x),\Exp^{-1}_\x(\y)\rangle_\x+\frac{C_f}{2}\rd^2(\x,\y).
$$
\end{definition}

Every point of a Riemannian manifold has a geodesically convex neighborhood.

\begin{definition}
A function $f:\cM\to\br$ is called \textbf{geodesically $L_f$-smooth} (\textbf{$L_f$-g-smooth})~\cite{zhang2016first-order} if $\|\cG f(\y)-\Gamma_\x^\y\cG f(\x)\|\le L_f \rd(\x,\y)$ for any $\x,\y\in\cM$.
Moreover, we have $f(\y)\leq f(\x)+\langle\cG f(\x),\Exp^{-1}_\x(\y)\rangle_\x+\frac{L_f}{2}\rd^2(\x,\y)$.
\end{definition}

On a Hadamard manifold the squared distance $\rd^2(\cdot,\y)$ is g-convex for every $\y$.

Under the LLS condition, the lower-level problem has a unique minimizer $\y^*(\x)$ for every $\x$, and the RBLO problem~\eqref{P:RBLO1} reduces to the single-level problem
\begin{equation}\label{P:RBLO2}
\min_{\x\in\cM}\vp(\x)=F(\x,\y^*(\x))\quad\text{s.t.}\quad\y^*(\x)=\arg\min_{\y\in\cN}f(\x,\y).
\end{equation}

\begin{definition}\label{D:hypergrad}
Under the LLS condition, the \textbf{hypergradient} at $\x\in\cM$ is the Riemannian gradient $\cG\vp(\x)\in\cT_\x\cM$ of the value function $\vp$ in~\eqref{P:RBLO2}, which accounts for the dependence of the lower-level solution $\y^*(\x)$ on $\x$~\cite{pedregosa2016hyperparameter}.
\end{definition}

The strong g-convexity also makes the lower Hessian $\Hess_\y f$ nonsingular, so $\y^*$ is differentiable in $\x$.
Differentiating the lower-level optimality condition $\cG_\y f(\x,\y^*(\x))=0$ along $\cM$ gives the Riemannian implicit function theorem
\begin{equation}\label{A:BRJ}
\D\y^*(\x)=-\bigl[\Hess_\y f(\x,\y^*(\x))\bigr]^{-1}\Hess_{\x\y}f(\x,\y^*(\x))\;:\;\cT_\x\cM\to\cT_{\y^*(\x)}\cN,
\end{equation}
and, substituting the adjoint of~\eqref{A:BRJ} into the chain rule $\cG\vp(\x)=\cG_\x F+[\D\y^*(\x)]^{*}\cG_\y F$, the ideal (implicit) hypergradient in closed form~\cite[Prop.~1]{han2024framework}, \cite[Prop.~10]{li2025riemannian}
\begin{equation}\label{A:HG1}
\cG\vp(\x)=\cG_\x F(\x,\y^*(\x))-\Hess_{\x\y}f(\x,\y^*(\x))^{*}\bigl[\Hess_\y f(\x,\y^*(\x))\bigr]^{-1}\cG_\y F(\x,\y^*(\x))\in\cT_\x\cM.
\end{equation}
The cross-Hessian $\Hess_{\x\y}f$ couples the two levels, while the inverse lower Hessian $[\Hess_\y f]^{-1}$ propagates the upper $\y$-gradient back to $\cT_\x\cM$.
Both are well defined precisely because the LLS condition keeps $\Hess_\y f$ invertible.

\begin{table}[t]
\centering\small
\renewcommand{\arraystretch}{1.2}
\begin{tabular}{@{}l l@{}}
\toprule
\textbf{Symbol} & \textbf{Meaning} \\
\midrule
$\cM,\ \cN$ & upper- and lower-level manifolds, connected and complete \\
$\x\in\cM,\ \y\in\cN$ & upper- and lower-level variables \\
$F,\ f:\cM\times\cN\to\br$ & upper- and lower-level objectives, both $C^2$ \\
$\cS(\x)$ & lower-level solution set, $\arg\min_{\y\in\cN}f(\x,\y)$ \\
$\y^*(\x)$ & the lower-level solution when $\cS(\x)$ is a singleton \\
$\vp(\x)$ & optimistic value, $\min_{\y\in\cS(\x)}F(\x,\y)$; the value function under LLS \\
$\Sopt(\x)$ & optimistic solution set, the minimizers of $F(\x,\cdot)$ over $\cS(\x)$ \\
$\Psi,\ \psi$ & $F(\x,\cdot)$ and $f(\x,\cdot)$ at a fixed $\x$ \\
$\cG f,\ \cG_\x f,\ \cG_\y f$ & Riemannian gradient and its two partial gradients \\
$\Hess_\y f,\ \Hess_{\x\y}f$ & lower-level Hessian and cross-Hessian \\
$\rd(\cdot,\cdot)$ & Riemannian distance \\
$\Exp_\y,\ \cR_\y$ & exponential map and retraction at $\y$ \\
\bottomrule
\end{tabular}
\caption{Symbols of the problem.}
\label{tab:symbols}
\end{table}

\section{Riemannian Bilevel Descent Aggregation}\label{sec:algorithm}

In this section, we present the proposed algorithm.
It applies the BDA mechanism as an aggregation strategy that aggregates the LL descent direction with UL information during the LL update.
The mechanism is adapted to the Riemannian setting through retraction-based updates.

\subsection{Unrolled hypergradient}\label{subsec:dynamics}
The hypergradient of Definition~\ref{D:hypergrad} requires the lower-level solution $\y^*(\x)$, and its closed form~\eqref{A:HG1} requires the inverse of the lower Hessian.
Neither is available when $\cS(\x)$ is not a singleton.
This happens when $f(\x,\cdot)$ is constant along a geodesic through a minimizer, for instance under a continuous symmetry of $f(\x,\cdot)$ or when $f$ depends on $\y$ only through a lower-dimensional quantity.
In this case, the value function is replaced by the optimistic value function $\vp(\x)=\inf_{\y\in\cS(\x)}F(\x,\y)$.

Since an explicit lower-level solution is out of reach in most cases, gradient-based bilevel methods~\cite{liu2022general,liu2021investigating} replace it by the iterate $\y_K(\x)$ reached after $K$ steps of an inner iteration and differentiate through these steps~\cite{mackay2019self,franceschi2017forward,franceschi2018bilevel,shaban2019truncated,nichol2018first}.
Unrolling the $K$ steps gives the problem
\begin{equation}\label{P:RBLO3}
    \min_{\x,\y_0,\dots,\y_K}\,F(\x,\y_K)\quad\text{s.t.}\quad \y_{k+1}=\cR_{\y_k}\bigl(-s_k\,g_k(\x)\bigr),\ k=0,\dots,K-1,
\end{equation}
where $s_k>0$ is a stepsize and $g_k(\x)\in\cT_{\y_k}\cN$ is the direction of the inner iteration.
Let $\vp_K(\x):=F(\x,\y_K(\x))$ denote the objective of~\eqref{P:RBLO3}.
The upper-level step uses the Riemannian gradient of $\vp_K$ as its descent direction,
\begin{equation}\label{A:HG2}
\widetilde{\cG}\vp(\x):=\cG\vp_K(\x)\in\cT_\x\cM,
\end{equation}
which we evaluate by reverse-mode automatic differentiation (AD)~\cite{griewank2008evaluating} through the $K$ unrolled updates, as in the reverse-mode hypergradient of~\cite{franceschi2017forward,franceschi2018bilevel}.

The direction $g_k$ determines which point of $\cS(\x)$ the iterate $\y_K(\x)$ approaches, and hence whether $\vp_K$ tracks the optimistic value $\vp$.
For a direction built from $\cG_\y f$ alone, the point of $\cS(\x)$ that the inner iterates approach depends only on the initialization.
\RBDA{} takes as inner iteration the aggregated update
\begin{equation}\label{A:agg}
    \y_{k+1}=\cR_{\y_k}\bigl(-s_k\,g_k\bigr),\qquad g_k=(1-\mu)\,\cG_\y f(\x,\y_k)+\mu\,\ak\,\cG_\y F(\x,\y_k),
\end{equation}
where $\x$ is fixed, $\mu\in[0,1]$ is the aggregation parameter~\cite{liu2022general}, $\ak>0$ is a multiplier with $\ak\to0$ as $k\to\infty$, and $\cR$ is a retraction on $\cN$.
For $\mu\in(0,1)$ the update~\eqref{A:agg} aggregates the LL and UL descent directions, so that the inner iterates are driven toward an optimistic solution in $\cS(\x)$.
Setting $\mu=0$ recovers the standard lower-level update, whereas $\mu=1$ ignores the LL objective and is not used.
Both gradients in~\eqref{A:agg} are evaluated at the same point $\y_k$ and therefore live in the same tangent space, so the aggregate is formed intrinsically and no transport is needed.
Since $f$ and $F$ are $C^2$ and $\cR$ is smooth, each step of~\eqref{A:agg} is $C^1$ in $(\x,\y_k)$, so for a fixed number of steps $K$ the objective $\vp_K$ is $C^1$ and the estimator~\eqref{A:HG2} is well defined.

Let $\cM=(-1,1)\subset\br$ and $\cN=\mathcal{S}^2\subset\br^3$ be the unit sphere, with $f(\x,\y)=\tfrac12(y_3-\x)^2$.
The minimizer set $\cS(\x)=\{\y\in\mathcal{S}^2:y_3=\x\}$ is an entire latitude circle, so LLS fails.
By rotational symmetry about the $z$-axis, $\cG_\y f$ has no component along this circle, so a purely lower-level solver keeps the longitude of $\y$ frozen at its initialization.
With $\mu\in(0,1)$, the upper-level gradient in~\eqref{A:agg} contributes a component along this circle and steers the iterates to an optimistic solution.
The sphere is positively curved and $f(\x,\cdot)$ is not g-convex on it, so this instance illustrates a non-singleton lower-level solution set but lies outside the Hadamard setting of the analysis.

The upper-level variable is updated by one Riemannian gradient step,
\begin{equation}\label{A:upd-upper}
\x^{r+1}=\cR_{\x^r}\!\big(-\eta\,\widetilde{\cG}\vp(\x^r)\big),
\end{equation}
with $r$ the outer index and $\eta>0$ the outer stepsize.
Algorithm~\ref{A:RBDA} summarizes the method with a general retraction.

\begin{algorithm}
\caption{Riemannian Bilevel Descent Aggregation (RBDA)}
\label{A:RBDA}
\begin{algorithmic}[1]
\State \textbf{Input:} upper-level initialization $\x^0$, lower-level initialization $\y_0$, aggregation parameter $\mu$, inner cap $K$, retraction $\cR$
\State \textbf{Input:} number of BB stepsize iterations $K_1\le K$, base stepsize $s$, decay exponent $\theta$, outer step $\eta$
\State $r\gets0$
\While{not converged}\Comment{outer stopping rule}
  \For{$k=0$ to $K-1$} \Comment{inner trajectory}
    \State $\ak\gets\tfrac12(k+1)^{-\theta}$ \Comment{multiplier~\eqref{A:ss2}}
    \State $g_k\gets(1-\mu)\,\cG_\y f(\x^r,\y_k)+\mu\,\ak\,\cG_\y F(\x^r,\y_k)$ \Comment{aggregated gradient~\eqref{A:agg}}
    \If{$1\le k\le K_1$}
      \State $s_k\gets\left|\langle\Delta\y_k,\Delta g_k\rangle_{\y_k}\right|/\langle\Delta g_k,\Delta g_k\rangle_{\y_k}$ \Comment{BB stepsize~\eqref{A:ss1}}
    \Else
      \State $s_k\gets s$
    \EndIf
    \State $\y_{k+1}\gets\cR_{\y_k}\!\bigl(-s_k\,g_k\bigr)$ \Comment{aggregated update~\eqref{A:agg}}
  \EndFor
  \State $\widetilde{\cG}\vp(\x^r)\gets$ reverse-mode AD of $F(\x^r,\y_K)$ through $\{\y_k\}_{k=0}^{K}$ \Comment{\eqref{A:HG2}}
  \State $\x^{r+1}\gets\cR_{\x^r}(-\eta\,\widetilde{\cG}\vp(\x^r))$;\quad $r\gets r+1$ \Comment{\eqref{A:upd-upper}}
\EndWhile
\State \textbf{Output:} $(\x^r,\y_K)$ of the last outer iteration
\end{algorithmic}
\end{algorithm}

\subsection{Stepsizes}\label{subsec:stepsize}
\RBDA{} runs the inner update~\eqref{A:agg} with
\begin{equation}\label{A:ss2}
    \ak=\tfrac12\,(k+1)^{-\theta},\qquad
    s_k=\begin{cases} s^{\mathrm{BB}}_k, & 1\le k\le K_1,\\ s, & \text{otherwise},\end{cases}
\end{equation}
where $s>0$ is a base stepsize and $\theta\in(\tfrac12,1]$ the decay exponent.
With~\eqref{A:ss2}, the update~\eqref{A:agg} is the aggregated scheme of~\cite[Eq.~(9)]{liu2022general} with a constant lower-level coefficient, whose analysis takes $\ak=1/(k+1)$, the case $\theta=1$.

On the iterations $1\le k\le K_1$ the stepsize is a nonmonotone Barzilai--Borwein (BB) stepsize~\cite{barzilai1988two,raydan1993barzilai}, computed from consecutive iterates and aggregated gradients,
\begin{equation}\label{A:ss1}
    s^{\mathrm{BB}}_k=\frac{\left|\bigl\langle\Delta\y_k,\Delta g_k\bigr\rangle_{\y_k}\right|}{\bigl\langle\Delta g_k,\Delta g_k\bigr\rangle_{\y_k}},\qquad
    \Delta \y_k=\Gamma_{\y_{k-1}}^{\y_k}\cR_{\y_{k-1}}^{-1}(\y_k),\qquad
    \Delta g_k=g_k-\Gamma_{\y_{k-1}}^{\y_k}\,g_{k-1},
\end{equation}
which is the second BB formula transported to the current tangent space.

\section{Convergence Analysis}\label{sec:analysis}
In this section, we analyze the inner iteration of \RBDA{} at a fixed upper-level variable $\x$ when the lower-level objective is g-convex on a Hadamard manifold and its solution set need not be a singleton.
Section~\ref{subsec:setting} states the assumptions, and Section~\ref{subsec:selection} gives the results.
Theorem~\ref{AN:main} shows that the inner iterates converge to a point of the optimistic solution set, and Theorem~\ref{AN:rate} gives rates under a growth condition of the upper-level objective on the lower-level solution set.
Theorem~\ref{AN:consist} shows that approximate minimizers of the objective with a finite number of inner iterations converge to minimizers of the optimistic value as this number grows.
Corollary~\ref{AN:lls} treats the LLS condition, under which the last inner iterate converges to the lower-level solution.
The results concern the inner loop behind each hypergradient, and they give no complexity bound for the upper-level iteration.
All proofs are given in Appendix~\ref{app:proofs}.

\subsection{Setting and assumptions}\label{subsec:setting}
Fix $\x\in\cM$ and write $\Psi:=F(\x,\cdot)$ and $\psi:=f(\x,\cdot)$, with $\cS:=\cS(\x)=\arg\min_{\y\in\cN}\psi$ and $\psi^*:=\min_{\cN}\psi$.
The function $\psi$ need not be strongly g-convex, so $\cS$ need not be a singleton.

With the exponential map as retraction and a constant stepsize $s$, the aggregated update~\eqref{A:agg} becomes
\begin{equation}\label{AN:upd}
\y_{k+1}=\Exp_{\y_k}\!\bigl(-s\,\bigl((1-\mu)\,\cG\psi(\y_k)+\mu\,\ak\,\cG\Psi(\y_k)\bigr)\bigr),
\end{equation}
where $\ak$ is the multiplier in~\eqref{A:ss2}.
Algorithm~\ref{A:RBDA} allows a general retraction $\cR$, whereas the analysis uses the exponential map.
Consequently, the results of this section apply to Algorithm~\ref{A:RBDA} with $\cR=\Exp$, and the case of a general retraction is not covered.
The analysis concerns the update~\eqref{AN:upd} with the constant stepsize $s$, and the first $K_1$ inner iterations of Algorithm~\ref{A:RBDA}, which use the BB stepsize, are treated at the end of Section~\ref{subsec:selection}.
We address three questions about the sequence $\{\y_k\}$ generated by~\eqref{AN:upd}.
The first is whether the iterates converge to a single point of $\Sopt(\x)$.
The second is at what rate $\psi(\y_k)$ approaches $\psi^*$ and $\Psi(\y_k)$ approaches $\vp(\x)$.
The third is whether, as $K\to\infty$, approximate minimizers of the objective $\vp_K$ of~\eqref{P:RBLO3} converge to minimizers of $\vp$.

To answer these questions, we impose conditions on the lower-level manifold $\cN$ and on the two objectives at the fixed $\x$.
The first assumption concerns the curvature of $\cN$, the convexity of $\psi$ and the existence of an optimistic solution.

\begin{assumption}\label{AN:geom}
\emph{(i)} $\cN$ is a Hadamard manifold whose sectional curvature is bounded below by a constant $\kappa^-\le0$.
\emph{(ii)} $\psi$ is g-convex on $\cN$ and $\cS$ is nonempty.
\emph{(iii)} $\Psi$ attains its infimum over $\cS$, so that $\Sopt(\x)$ is nonempty.
\end{assumption}

Since a Hadamard manifold is noncompact, Assumption~\ref{AN:geom} excludes compact lower-level manifolds.
The assumption concerns the lower-level manifold alone, and the upper-level manifold $\cM$ may be compact.
Convexity of $\psi$ is required on all of $\cN$ because $\cS$ and $\psi^*$ are global objects.
The set $\cS$ is closed and g-convex as a sublevel set of the continuous g-convex function $\psi$.
Hence the metric projection $P_\cS(\y):=\arg\min_{\y'\in\cS}\rd(\y,\y')$ is single-valued and nonexpansive~\cite[Prop.~II.2.4]{bridson1999metric}.
Condition \emph{(iii)} holds, for instance, when $\cS$ is bounded, since a closed and bounded subset of a complete Riemannian manifold is compact~\cite[Ch.~6]{lee2018introduction}.

By Assumption~\ref{AN:geom}\emph{(iii)} we may fix an optimistic solution $\by\in\Sopt(\x)$ and a radius $\rho>0$, and define the closed geodesic balls
\begin{equation}\label{AN:ball}
\bB:=\{\y\in\cN:\rd(\y,\by)\le\rho\},\qquad \bB_2:=\{\y\in\cN:\rd(\y,\by)\le2\rho\}.
\end{equation}
On a Hadamard manifold both balls are compact and g-convex, and the local constants of the analysis are defined on them.

The analysis of the update~\eqref{AN:upd} uses a comparison inequality for the squared distance.
In Euclidean space, where $\Exp_\y(v)=\y+v$, the expansion $\|\y+v-\y'\|^2=\|\y-\y'\|^2-2\langle v,\y'-\y\rangle+\|v\|^2$ is an identity.
On $\cN$ the same expansion holds as an inequality in which the squared step $\|v\|^2$ is multiplied by a factor depending on the curvature.

\begin{lemma}[{\cite[Lemma~5]{zhang2016first-order}}]\label{AN:cmplemma}
Let Assumption~\ref{AN:geom} hold.
For all $\y,\y'\in\bB$ and every $v\in\cT_\y\cN$,
\begin{equation}\label{AN:cmp}
\rd^2\big(\Exp_{\y}(v),\,\y'\big)\;\le\;\rd^2(\y,\y')-2\big\langle v,\Exp^{-1}_{\y}(\y')\big\rangle_{\y}+\zeta\,\|v\|_{\y}^2,
\end{equation}
where $\zeta:=2\rho\sqrt{|\kappa^-|}/\tanh(2\rho\sqrt{|\kappa^-|})$ for $\kappa^-<0$ and $\zeta:=1$ for $\kappa^-=0$.
\end{lemma}

The cited lemma gives~\eqref{AN:cmp} with $2\rho$ replaced by $\rd(\y,\y')$, and its factor is nondecreasing in this distance, which is at most $2\rho$ on $\bB$.
The factor $\zeta\ge1$ multiplies only the squared step $\|v\|^2$, whereas the coefficient of the distance term $\rd^2(\y,\y')$ is one.
Inequality~\eqref{AN:cmp} holds for every $v\in\cT_\y\cN$, including those for which $\Exp_\y(v)$ lies outside $\bB$.
This fact is used to prove that the iterates remain in $\bB$.

\begin{assumption}\label{AN:smooth}
\emph{(i)} On $\bB_2$ the function $\psi$ is $L_f$-g-smooth, with $L_f\ge G_f/\rho$ and $G_f:=\max_{\bB}\|\cG\psi\|$.
\emph{(ii)} $\Psi$ is g-convex on $\bB$, and we write $G_F:=\max_{\bB}\|\cG\Psi\|$.
\end{assumption}

The condition $L_f\ge G_f/\rho$ restricts nothing, since a g-smoothness constant may be increased, and it ensures that the point $\Exp_\y(-\cG\psi(\y)/L_f)$ lies in $\bB_2$ for every $\y\in\bB$.

Since g-convexity of $\psi$ gives no bound on $\rd(\y,\cS)$ in terms of $\psi(\y)-\psi^*$, we require the lower level to satisfy a quadratic growth condition.

\begin{assumption}\label{AN:qg}
There is a constant $c>0$ such that
\begin{equation}\label{AN:qgineq}
\rd^2(\y,\cS)\le c\,\big(\psi(\y)-\psi^*\big)\qquad\text{for every }\y\in\bB.
\end{equation}
\end{assumption}

Inequality~\eqref{AN:qgineq} is the quadratic growth condition of $\psi$ relative to its solution set~\cite[Cor.~3.6]{drusvyatskiy2018error}, and its Riemannian form is given in~\cite[Def.~1.2]{rebjock2024fast}.
For $C^2$ functions on a Riemannian manifold, quadratic growth and the Polyak--\L{}ojasiewicz condition are equivalent near each point of $\cS$, up to an arbitrarily small decrease of the constant~\cite[Table~1]{rebjock2024fast}.
Quadratic growth of the lower-level objective is also assumed in the penalty-based bilevel gradient descent method of~\cite[Assumption~4(i)]{shen2023penalty}.
If $\psi$ is $C_f$-strongly g-convex, its solution set is $\cS=\{\y^*(\x)\}$, and the strong g-convexity inequality of Section~\ref{sec:preliminaries} applied at $\y^*(\x)$, where $\cG\psi$ vanishes, gives~\eqref{AN:qgineq} with $c=2/C_f$.

\begin{assumption}\label{AN:radius}
The initial point satisfies $\rd(\y_0,\by)<\rho$.
\end{assumption}

\subsection{Convergence}\label{subsec:selection}
The results of this section use the following conditions on the parameters of~\eqref{AN:upd}:
\begin{subequations}\label{AN:param}
\begin{align}
&\mu\in(0,1),\qquad 0<s\le\frac{1}{4\zeta(1-\mu)L_f},\qquad \rd^2(\y_0,\by)+C_0\sum_{k\ge0}\ak^2\le\rho^2,\label{AN:parama}\\
&\{\ak\}\subset(0,1],\qquad \ak\to0,\qquad \sum_{k\ge0}\ak=\infty,\qquad \sum_{k\ge0}\ak^2<\infty,\label{AN:paramb}
\end{align}
\end{subequations}
where $C_0:=2\zeta s^2\mu^2G_F^2+2s\mu^2G_F^2c/(1-\mu)$.
The following theorem answers the first question, and it shows that the inner iterates converge to a point of $\Sopt(\x)$.

\begin{theorem}\label{AN:main}
Let Assumptions~\ref{AN:geom}--\ref{AN:radius} hold, and let $\mu$, $\{\ak\}$ and $s$ satisfy~\eqref{AN:param}.
Then the iterates of~\eqref{AN:upd} satisfy the following:
\begin{enumerate}[label=\textup{(\roman*)}]
  \item $\lim_{k\to\infty}\psi(\y_k)=\psi^*$ and $\lim_{k\to\infty}\Psi(\y_k)=\vp(\x)$;
  \item $\lim_{k\to\infty}\y_k=\hat\y$ for some $\hat\y\in\Sopt(\x)\cap\bB$.
\end{enumerate}
\end{theorem}

Its proof is deferred to Appendix~\ref{app:pfmain}.
The multipliers~\eqref{A:ss2} satisfy~\eqref{AN:paramb} for every $\theta\in(\tfrac12,1]$.
Since $C_0$ tends to $0$ as $s\to0$, Assumption~\ref{AN:radius} ensures that the last inequality in~\eqref{AN:parama} holds for every sufficiently small $s$.
Assumption~\ref{AN:qg} is used in the proof to bound the amount by which $\Psi(\y_k)$ falls below $\vp(\x)$ while $\y_k$ lies outside $\cS$.

Theorem~\ref{AN:main} gives no rate, and the second question needs more information about $\Psi$ on $\cS$.
Since $\Psi$ may increase arbitrarily slowly on $\cS$ away from $\Sopt(\x)$, the distance to $\Sopt(\x)$ admits a rate only if this increase is quantified.
The following theorem quantifies this increase by a growth condition of order $p$ and gives rates for the multipliers~\eqref{A:ss2} with $\theta<1$.

\begin{theorem}\label{AN:rate}
Let Assumptions~\ref{AN:geom}--\ref{AN:radius} hold, let $a(k+1)^{-\theta}\le\ak\le(k+1)^{-\theta}$ for every $k$, where $a\in(0,1]$ and $\theta\in(\tfrac12,1)$, and let $\mu$ and $s$ satisfy~\eqref{AN:parama}.
Suppose that there are constants $c'>0$ and $p\ge1$ such that $\Psi(\y)-\vp(\x)\ge c'\,\rd^{\,p}(\y,\Sopt(\x))$ for every $\y\in\cS\cap\bB$.
Then for every $k\ge2$ there is an index $j_k$ with $k/2\le j_k\le k$ such that the following hold:
\begin{enumerate}[label=\textup{(\roman*)}]
  \item $\psi(\y_{j_k})-\psi^*=O(k^{-1})$ and $\rd(\y_{j_k},\cS)=O(k^{-1/2})$;
  \item $|\Psi(\y_{j_k})-\vp(\x)|=O(k^{-(1-\theta)})$ and $\rd(\y_{j_k},\Sopt(\x))=O(k^{-(1-\theta)/p})$,
\end{enumerate}
where the constants in $O(\cdot)$ depend on $\rho$, $s$, $\mu$, $a$, $c$, $c'$, $p$ and $G_F$, and are independent of $k$.
\end{theorem}

Its proof is deferred to Appendix~\ref{app:pfrate}, where the constants are given explicitly.
The estimates hold at an index $j_k$ in the window $[k/2,k]$, which need not be $k$.
Estimate~\emph{(i)} does not use the growth condition on $\Psi$.
For $p=1$ that condition is a weak sharp minimum of $\Psi$ over $\cS$ in the sense of~\cite[eq.~(1)]{burke1993weak}, and for general $p$ it is the H\"olderian growth of~\cite[Sec.~2.4]{bolte2017error}.
In~\emph{(ii)}, a smaller $\theta$ gives a larger exponent $1-\theta$, whereas the bound $\sum_{k\ge0}(k+1)^{-2\theta}$ on $\sum_{k\ge0}\ak^2$ grows as $\theta\downarrow\tfrac12$, so that~\eqref{AN:parama} requires a smaller $s$.

The growth condition on $\Psi$ holds, for instance, when $\Psi$ is strongly g-convex on $\bB$.

\begin{proposition}\label{AN:sgrowth}
Let Assumption~\ref{AN:geom} hold, and let $\Psi$ be $C_F$-strongly g-convex on $\bB$.
Then $\Sopt(\x)\cap\bB=\{\by\}$, and
$$
\Psi(\y)-\vp(\x)\ge\frac{C_F}{2}\,\rd^2\big(\y,\Sopt(\x)\big)\qquad\text{for every }\y\in\cS\cap\bB,
$$
that is, the growth condition of Theorem~\ref{AN:rate} holds with $p=2$ and $c'=C_F/2$.
\end{proposition}

Its proof is deferred to Appendix~\ref{app:pfsgrowth}.
Under the assumptions of Proposition~\ref{AN:sgrowth}, the limit in Theorem~\ref{AN:main}\emph{(ii)} is $\by$, and Theorem~\ref{AN:rate}\emph{(ii)} gives $\rd(\y_{j_k},\by)=O(k^{-(1-\theta)/2})$.

At a fixed $\x$, the inner problem is the simple bilevel problem of minimizing $\Psi$ over the solution set of $\psi$.
For convex simple bilevel problems in Euclidean space, the bi-sub-gradient method of~\cite[Theorem~4.11]{merchav2023convex} attains the rates $O(k^{-\beta})$ and $O(k^{-(1-\beta)})$ for the lower- and upper-level objective values with $\beta\in(\tfrac12,1)$, at an index in the window $[k,2k]$.
The upper-level estimate in Theorem~\ref{AN:rate} has the same form $O(k^{-(1-\theta)})$, whereas the lower-level estimate has the exponent $1$ for every $\theta$, and Theorem~\ref{AN:rate} gives in addition rates for $\rd(\y_{j_k},\cS)$ and $\rd(\y_{j_k},\Sopt(\x))$.
These estimates use Assumption~\ref{AN:qg} and the growth condition on $\Psi$, which~\cite{merchav2023convex} does not require.
In both results the index is only shown to exist, and identifying $j_k$ requires $\psi^*$ and $\vp(\x)$.
For the Euclidean aggregation scheme with $\ak=1/(k+1)$, that is, $\theta=1$, \cite[Theorem~4]{liu2022general} gives a rate of order $\big((1+\ln k)/k^{1/4}\big)^{1/2}$ for the lower-level objective and no rate for the upper-level objective.

Theorems~\ref{AN:main} and~\ref{AN:rate} concern a fixed upper-level variable, whereas the third question concerns the dependence on $\x$.
For $\x\in\cM$, let $\y_K(\x)$ denote the $K$-th iterate of~\eqref{AN:upd} at $\x$ from a starting point $\y_0\in\cN$ that does not depend on $\x$, and let $\vp_K(\x):=F(\x,\y_K(\x))$ as in~\eqref{P:RBLO3}.
The following theorem shows that approximate minimizers of $\vp_K$ over a compact set converge to minimizers of $\vp$ when Assumptions~\ref{AN:geom}--\ref{AN:radius} and~\eqref{AN:param} hold with constants independent of $\x$ on that set.

\begin{theorem}\label{AN:consist}
Let $\cX\subset\cM$ be nonempty and compact.
Suppose that for every $\x\in\cX$, Assumptions~\ref{AN:geom}--\ref{AN:radius} and~\eqref{AN:param} hold with the same $\rho$, $L_f$, $G_F$, $c$, $\mu$, $s$ and $\{\ak\}$, where the center $\by(\x)\in\Sopt(\x)$ may depend on $\x$.
Let $\x_K\in\cX$ satisfy $\vp_K(\x_K)\le\vp_K(\x)+\varepsilon_K$ for every $\x\in\cX$, where $\varepsilon_K\to0$.
Then the following hold:
\begin{enumerate}[label=\textup{(\roman*)}]
  \item every limit point of $\{\x_K\}$ minimizes $\vp$ over $\cX$;
  \item $\lim_{K\to\infty}\inf_{\x\in\cX}\vp_K(\x)=\inf_{\x\in\cX}\vp(\x)$.
\end{enumerate}
\end{theorem}

Its proof is deferred to Appendix~\ref{app:pfconsist}.
Theorem~\ref{AN:consist} is the Riemannian counterpart of~\cite[Theorem~1]{liu2022general}, which extends~\cite[Theorem~1]{liu2020generic} from exact minimizers to $\varepsilon_K$-minimizers.
The two hypotheses of~\cite[Theorem~1]{liu2022general}, the uniform convergence of the lower-level gap on $\cX$ with bounded iterates and the pointwise convergence of $\vp_K$ to $\vp$, follow here from the uniformity of the constants in Theorem~\ref{AN:main}.
The proof uses a bound on the lower-level gap at the last iterate $\y_K(\x)$, and the independence of the constants from $\x$ makes this bound uniform on $\cX$.
Theorem~\ref{AN:consist} concerns approximate minimizers of $\vp_K$ over $\cX$, and it gives no estimate for the iterates of the upper-level update~\eqref{A:upd-upper}.

Under the LLS condition, $\cS=\{\y^*(\x)\}$, and the following corollary gives a rate for the last iterate.
In this case Assumption~\ref{AN:qg} holds with $c=2/C_f$, and the g-convexity of $\Psi$ in Assumption~\ref{AN:smooth}\emph{(ii)} is not required.

\begin{corollary}\label{AN:lls}
Let Assumption~\ref{AN:geom}\emph{(i)} hold, and let $\psi$ be g-convex on $\cN$ with a unique minimizer $\y^*(\x)$.
Let the balls~\eqref{AN:ball} be centered at $\by=\y^*(\x)$, let $\psi$ be $C_f$-strongly g-convex on $\bB_2$ and satisfy Assumption~\ref{AN:smooth}\emph{(i)}, let $G_F:=\max_{\bB}\|\cG\Psi\|$, and let Assumption~\ref{AN:radius} hold.
Let $\mu$, $\{\ak\}$ and $s$ satisfy~\eqref{AN:param} with $c=2/C_f$, and let $q:=1-s(1-\mu)C_f/4$.
Then $q\in[\tfrac{15}{16},1)$, and the following hold:
\begin{enumerate}[label=\textup{(\roman*)}]
  \item $\rd^2(\y_k,\y^*(\x))\le q^{k/2-1}\rho^2+\dfrac{C_0}{1-q}\max_{j\ge k/2}\alpha_j^2$ for every $k\ge2$;
  \item if $\ak\le a(k+1)^{-\theta}$ for every $k$, where $a\in(0,1]$ and $\theta\in(\tfrac12,1]$, then $\rd(\y_k,\y^*(\x))=O(k^{-\theta})$,
\end{enumerate}
where the constant in $O(\cdot)$ depends on $\rho$, $s$, $\mu$, $a$, $\theta$, $\zeta$, $G_F$ and $C_f$, and is independent of $k$.
\end{corollary}

Its proof is deferred to Appendix~\ref{app:pflls}.
Unlike Theorem~\ref{AN:rate}, Corollary~\ref{AN:lls} bounds the last iterate and includes the multipliers~\eqref{A:ss2} with $\theta=1$.
The first term in~\emph{(i)} decreases linearly, whereas the second is of order $\max_{j\ge k/2}\alpha_j^2$, since the upper-level term $s\mu\ak\cG\Psi(\y_k)$ in~\eqref{AN:upd} does not vanish at $\y^*(\x)$.
Consequently, the exponent in~\emph{(ii)} is that of the multiplier.
For $\mu=0$ in~\eqref{AN:upd}, the recursion in Appendix~\ref{app:pflls} has no second term and $\rd^2(\y_k,\y^*(\x))$ decreases linearly, so under the LLS condition the aggregation replaces a linear rate by the rate $O(k^{-\theta})$ in~\emph{(ii)}.
The curvature of $\cN$ enters the analysis only through the factor $\zeta$ in the bound on $s$ and in $C_0$, and the exponents in Theorem~\ref{AN:rate} and in~\emph{(ii)} do not depend on it.

Algorithm~\ref{A:RBDA} uses the BB stepsize on the inner iterations $1\le k\le K_1$ and the update~\eqref{AN:upd} afterwards.
For $K_1\ge1$ and the multipliers~\eqref{A:ss2}, the sequence $\ty_j:=\y_{K_1+1+j}$ is generated by~\eqref{AN:upd} from $\ty_0=\y_{K_1+1}$ with the multipliers $\tilde\alpha_j:=\alpha_{K_1+1+j}=\tfrac12(K_1+2+j)^{-\theta}$.
Since $j+1\le K_1+2+j\le(K_1+2)(j+1)$, these multipliers satisfy $a_1(j+1)^{-\theta}\le\tilde\alpha_j\le\tfrac12(j+1)^{-\theta}$ with $a_1:=\tfrac12(K_1+2)^{-\theta}$.
Hence, if Assumptions~\ref{AN:geom}--\ref{AN:radius} and~\eqref{AN:param} hold with $\ty_0$ and $\{\tilde\alpha_j\}$ in place of $\y_0$ and $\{\ak\}$, Theorem~\ref{AN:main} and Corollary~\ref{AN:lls} apply to $\{\ty_j\}$, and Theorem~\ref{AN:rate} applies to $\{\ty_j\}$ with $a=a_1$ when $\theta<1$.

The inequality $i^{-2\theta}\le\int_{i-1}^{i}t^{-2\theta}\,dt$ gives
$$
\sum_{j\ge0}\tilde\alpha_j^2=\frac14\sum_{i\ge K_1+2}i^{-2\theta}\le\frac14\int_{K_1+1}^{\infty}t^{-2\theta}\,dt=\frac{(K_1+1)^{1-2\theta}}{4(2\theta-1)} ,
$$
so the sum in~\eqref{AN:parama} decreases as $K_1$ grows, and a larger $K_1$ admits a larger stepsize $s$, whereas the coefficient $a_1$ in the constants of Theorem~\ref{AN:rate} decreases.
Theorem~\ref{AN:consist} is stated for $K_1=0$, since its proof uses a starting point that does not depend on $\x$, whereas $\ty_0=\y_{K_1+1}(\x)$ depends on $\x$ for $K_1\ge1$.

\section{Numerical Experiments}\label{sec:experiments}
We report four experiments, the first two synthetic and the last two on real data.
Section~\ref{subsec:lls} considers a lower-level singleton problem, on which every estimator is valid, so that the comparison concerns cost alone.
Section~\ref{subsec:nonlls} adds to it a lower-level block that the lower-level objective ignores, so that the lower-level solution set is no longer a singleton.
Section~\ref{subsec:fewshot} considers few-shot meta-learning, whose one-shot episodes underdetermine the lower-level heads.
Section~\ref{subsec:real} considers data hyper-cleaning with no lower-level regularization on three MNIST-family datasets.
Appendix~\ref{app:supp} collects the supplementary experiments.

We compare \RBDA{} with three hypergradient estimators built on the plain lower-level update and with two recent RBLO solvers.
CG and NS solve the implicit system of the hypergradient~\eqref{A:HG1} by $T$ Hessian-vector iterations of conjugate gradients and of the truncated Neumann series, respectively.
AD differentiates through the unrolled trajectory of the plain lower-level update by reverse-mode automatic differentiation.
\RBDA{} is Algorithm~\ref{A:RBDA}.
AD and \RBDA{} store the inner trajectory, so their memory grows linearly in the inner cap, whereas CG and NS need $O(1)$ storage.
RF$^2$SA~\cite{dutta2024riemannian} is fully first-order and runs the multiplier schedule $\lambda_k=\lambda_0(1+k)^{1/3}$ of~\cite[Cor.~3]{dutta2024riemannian}.
AdaRHD~\cite{shi2025adaptive} restarts its inner accumulator at each outer step and uses the adaptive outer stepsize of its Algorithm~1.
The outer stepsize of RF$^2$SA is the one defined by its multiplier schedule, so the fixed outer stepsize $\eta_x$ below applies to CG, NS, AD and \RBDA{} only.
Both RF$^2$SA and AdaRHD presume a strongly g-convex lower level, which the problem of Section~\ref{subsec:lls} satisfies and the other three problems do not.

We write $R$ for the number of outer steps and $K$ for the inner cap of Algorithm~\ref{A:RBDA}.
Here $\eta_x$ is the outer stepsize $\eta$ of Algorithm~\ref{A:RBDA}, $\eta_y$ its base stepsize $s$, and $T$ the truncation length of CG and NS.
Every inner loop starts from the terminal iterate of the previous outer step.
The constants of Algorithm~\ref{A:RBDA} are the same on the four problems: $\theta=\tfrac34$, $\mu=0.7$, $K_1=K/2$, and the BB stepsize is clipped to $[10^{-4},s_{\max}]$ and treated as a constant by the reverse sweep.
The comparisons run on $20$ random seeds and report mean $\pm$ standard deviation over seeds.

The implementation is in \texttt{PyTorch}, with \texttt{geoopt}~\cite{kochurov2020geoopt} providing the manifold primitives.
All runs are computed on CPU on Google Colab.
The code and the result files are available at \url{https://github.com/Cz1544252489/RBDA-Public}.

\subsection{Lower-level singleton baseline}\label{subsec:lls}
This experiment measures the cost of the aggregation when the lower-level solution is unique.
We consider the synthetic problem of~\cite{han2024framework},
\begin{equation}\label{P:syn}
\begin{gathered}
\max_{W\in\St(d,r)}\ \tr\big(M^\star(W)\,X^\top YW^\top\big)\\
\text{s.t.}\quad
M^\star(W)=\arg\min_{M\in\cS_{++}^d}\ \langle M,X^\top X\rangle+\langle M^{-1},WY^\top YW^\top+\nu I\rangle,
\end{gathered}
\end{equation}
with data $X\in\br^{n\times d}$, $Y\in\br^{n\times r}$ and $\nu>0$.
The upper-level manifold is the Stiefel manifold $\St(d,r)=\{W\in\br^{d\times r}:W^\top W=I_r\}$, which is compact.
The lower-level manifold is the set $\cS_{++}^d$ of symmetric positive definite $d\times d$ matrices with the affine-invariant metric $\langle U,V\rangle_M=\tr(M^{-1}UM^{-1}V)$, a Hadamard manifold~\cite{sra2015conic} whose sectional curvature is bounded below by homogeneity.
Writing $A:=X^\top X$ and $B(W):=WY^\top YW^\top+\nu I$, the lower-level objective $\psi(M)=\langle M,A\rangle+\langle M^{-1},B(W)\rangle$ has the unique minimizer
$$
M^\star(W)=A^{-1/2}\big(A^{1/2}B(W)A^{1/2}\big)^{1/2}A^{-1/2},
$$
which supplies the reference lower-level solution and the exact hypergradient by~\eqref{A:HG1}.
Both terms of $\psi$ are g-convex on $\cS_{++}^d$~\cite[eq.~(2.3) and Ex.~2.5]{sra2015conic}, and along every geodesic the second derivative of $\psi$ is at least $\nu/\lambda_{\max}(M)$ times the squared speed, so $\psi$ is g-convex on $\cS_{++}^d$ with the unique minimizer $M^\star(W)$ and strongly g-convex on every bounded subset, as Corollary~\ref{AN:lls} requires.
The upper-level objective need not be g-convex in $M$, which Corollary~\ref{AN:lls} does not require.

We use $n=200$, $d=50$, $r=20$ and $\nu=0.01$, with the entries of $X$ and $Y$ standard normal divided by $\sqrt n$, $W_0$ uniform on $\St(d,r)$ and $M_0=I$.
We take $R=400$, $K=80$, $\eta_x=0.05$, $\eta_y=0.5$, $s_{\max}=10$ and $\lambda_0=1$.
The lower-level update of CG, NS, AD and \RBDA{} uses the retraction $\cR_M(v)=M+v+\tfrac12vM^{-1}v$ of \texttt{geoopt}, whereas AdaRHD and RF$^2$SA use the exponential map, and every method takes its upper-level step with the Cayley retraction on $\St(d,r)$~\cite{wen2013feasible}.
The analysis of Section~\ref{sec:analysis} covers the exponential map, and the runs of \RBDA{} in Table~\ref{tab:lls} use the retraction.
The inner loop stops once the lower-level gradient norm, or for \RBDA{} the length of the first update, falls below $10^{-3}$ of its value at the first iteration of that loop, and RF$^2$SA runs the full inner cap.
Every run is scored against a per-seed reference $U^\star$, obtained by ascent along the exact hypergradient from the same $W_0$ until the Riemannian gradient norm falls below $10^{-10}$.
The cost of a method is its solver time, the wall-clock time of the hypergradient evaluations and outer retractions, until the relative gap $(U^\star-U_r)/|U^\star|$ first falls below $\varepsilon=10^{-2}$ or $10^{-3}$.
The truncation length is $T=20$ for CG and $T=200$ for NS, since NS does not reach the $10^{-3}$ threshold at $T=20$.

\begin{table}[ht]
\centering
\renewcommand{\arraystretch}{1.12}
\small
\setlength{\tabcolsep}{4pt}
\begin{tabular}{@{}lcccc@{}}
\toprule
\textbf{Method} & \textbf{Time to $\varepsilon=10^{-2}$ (s)} & \textbf{Time to $\varepsilon=10^{-3}$ (s)} & \shortstack{\textbf{Outer steps}\\$10^{-2}$ / $10^{-3}$} & \shortstack{\textbf{Inner iterations}\\per outer step} \\
\midrule
CG ($T=20$) & $10.82\pm1.18$ & $21.25\pm2.07$ & $164$ / $311$ & $30.2$ \\
NS ($T=200$) & $22.81\pm1.98$ & $43.95\pm3.53$ & $166$ / $316$ & $30.1$ \\
AD & $13.91\pm1.43$ & $25.77\pm1.96$ & $164$ / $293$ & $30.2$ \\
RF$^2$SA~\cite{dutta2024riemannian} & --- & --- & --- / --- & $80$ \\
AdaRHD~\cite{shi2025adaptive} & $5.79\pm1.34$ & $11.01\pm1.58$ & $56$ / $111$ & $32.6$ \\
\RBDA{} & $40.11\pm4.62$ & $73.18\pm7.55$ & $162$ / $293$ & $47.7$ \\
\bottomrule
\end{tabular}
\caption{Singleton problem, $20$ seeds, $R=400$: solver time (mean $\pm$ standard deviation) until the relative gap to $U^\star$ first falls below $\varepsilon$, mean outer steps to each $\varepsilon$, and mean inner iterations per outer step.
Lower is better in the first four columns, and the last column is the cost of one outer step.
RF$^2$SA reaches neither threshold within $R=400$ on any seed.}
\label{tab:lls}
\end{table}

\begin{figure}[ht]
\centering
\includegraphics[width=\textwidth]{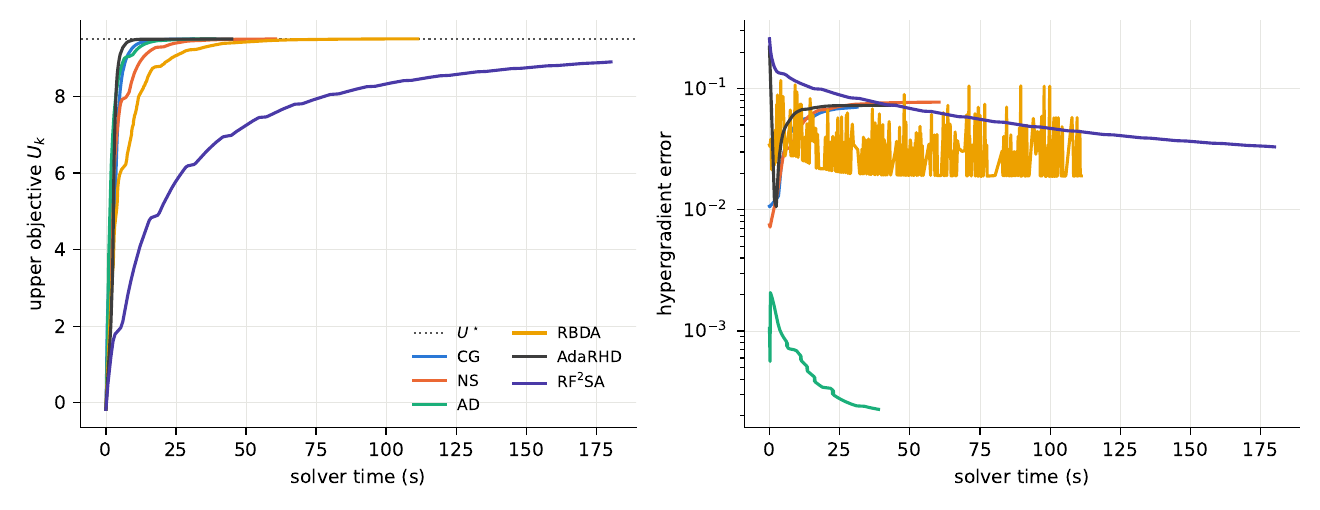}
\caption{Singleton problem, first seed, $R=400$, against solver time: (a) upper-level objective with the reference $U^\star$ dotted, and (b) hypergradient error, the norm of the difference between the estimated and the exact hypergradient.}
\label{fig:lls}
\end{figure}

The four fixed-stepsize estimators follow the same outer trajectory, as Table~\ref{tab:lls} and Figure~\ref{fig:lls} show, and differ only in the cost per outer step.
NS needs twice the time of CG because of its ten times longer truncation, and the extra time of AD is the reverse sweep through the unrolled loop.
\RBDA{} costs about three times AD, since with $\mu=0.7$ the aggregated inner loop takes more iterations than the plain lower-level update to satisfy the same stopping rule and each iteration evaluates both gradients.
Its hypergradient carries the bias of the aggregation, two orders of magnitude above the error of AD, which does not slow the outer iteration at these thresholds.
AdaRHD is the fastest, since its adaptive outer step reaches both thresholds in about a third of the outer steps.
RF$^2$SA reaches neither threshold, since its outer stepsize decreases as the multiplier $\lambda_k$ grows.
The BB stepsize saves about a fifth of the time to the $10^{-3}$ threshold, and Appendix~\ref{app:bb} reports this ablation.

Table~\ref{tab:geometry} isolates the lower-level problem of~\eqref{P:syn} at the fixed upper variable $W_0$ with $\mu=0$, with the columns of $X$ scaled so that $\operatorname{cond}(A)$ ranges from $8.5$ to $1.4\times10^3$.
E-BDA is the lower-level descent of BDA~\cite{liu2022general} in the Euclidean geometry, the Euclidean gradient with the additive update $M+v$.
\RBDA{}-retr uses the affine-invariant gradient with the same additive update, and \RBDA{} the affine-invariant gradient with the exponential map.
Each update uses the stepsize of a grid that reaches $\rd(M,M^\star(W_0))\le10^{-2}$ in the fewest steps.
E-BDA reaches the target at $\kappa=1$ and $\kappa=10$ only, whereas the two affine-invariant updates reach it in tens of steps at every $\kappa$ and within a factor $1.3$ of each other.
The metric acts as a preconditioner of the gradient iteration~\cite{mishra2016riemannian}, and across the four values of $\kappa$ the condition number of the lower-level Hessian at $M^\star(W_0)$ grows by a factor of about $500$ in the Euclidean metric and about $11$ in the affine-invariant metric.

\begin{table}[ht]
\centering
\renewcommand{\arraystretch}{1.12}
\small
\setlength{\tabcolsep}{5pt}
\begin{tabular}{@{}lcccccc@{}}
\toprule
& & \multicolumn{2}{c}{\textbf{Condition number}} & \multicolumn{3}{c}{\textbf{Inner steps to $\rd(M,M^\star)\le10^{-2}$}} \\
\cmidrule(lr){3-4}\cmidrule(lr){5-7}
$\kappa$ & $\operatorname{cond}(A)$ & Euclidean & affine-invariant & E-BDA & \RBDA{}-retr & \RBDA{} \\
\midrule
$1$ & $8.5$ & $70$ & $24$ & $219$ & $23$ & $20$ \\
$10$ & $25.6$ & $243$ & $41$ & $1{,}409$ & $37$ & $29$ \\
$10^2$ & $172$ & $2.3\times10^3$ & $99$ & --- & $60$ & $50$ \\
$10^3$ & $1{,}440$ & $3.6\times10^4$ & $264$ & --- & $78$ & $79$ \\
\bottomrule
\end{tabular}
\caption{Lower-level geometry at the fixed upper variable $W_0$ on the first seed with $\mu=0$ and the BB stepsize: condition number of the lower-level Hessian at $M^\star(W_0)$ in the two metrics, and inner steps of the three updates from $M_0=I$ to $\rd(M,M^\star(W_0))\le10^{-2}$ at their best stepsize in $\{1,0.5,0.2,0.1,\dots,10^{-3}\}$.
A dash means that no stepsize of the grid reaches the target within $2{,}000$ steps.}
\label{tab:geometry}
\end{table}

\subsection{Ignored-block counterexample}\label{subsec:nonlls}
This experiment tests the selection when the lower-level solution set is a positive-dimensional submanifold.
We add to the lower-level variable of~\eqref{P:syn} a second block $P\in\cS_{++}^d$ on which the lower-level objective does not depend, and minimize
\begin{equation}\label{P:nonlls}
\begin{gathered}
\min_{W\in\St(d,r)}\ \vp(W)=\min_{(M,P)\in\cS(W)}\ \langle M,G(W)\rangle+\tr(P)+\langle P^{-1},\Sigma(W)\rangle\\
\text{s.t.}\quad
\cS(W)=\arg\min_{M,P\in\cS_{++}^d}\ \langle M,X^\top X\rangle+\langle M^{-1},WY^\top YW^\top+\nu I\rangle,
\end{gathered}
\end{equation}
with $G(W)=\operatorname{sym}(X^\top YW^\top)$ and $\Sigma(W)=WY^\top YW^\top+\omega I$.
This is the Riemannian counterpart of the ignored-block counterexamples of~\cite[Ex.~1, Rem.~3]{liu2022general}.
The lower-level solution set is $\cS(W)=\{M^\star(W)\}\times\cS_{++}^d$, and the upper-level objective selects $P^\star(W)=\Sigma(W)^{1/2}$ in it, so that $\vp(W)$, $\rd(P,P^\star(W))$ and $\rd(M,M^\star(W))$ are available in closed form.
Since the lower-level gradient has no component along $P$, every inner loop driven by it keeps $P$ at its initial value, and CG, NS, AD and AdaRHD stay at the initial selection gap $\rd(P_0,P^\star(W))$ whatever their accuracy in $M$.
The problem satisfies Assumptions~\ref{AN:geom}, \ref{AN:qg} and~\ref{AN:radius} at every $W$, whereas the coupling $\langle M,G(W)\rangle$ has an indefinite coefficient and the upper-level objective is not g-convex in $M$, so Assumption~\ref{AN:smooth}\emph{(ii)} fails and Theorem~\ref{AN:main} does not cover it.

We use $n=80$, $d=10$, $r=5$, $\nu=0.01$ and $\omega=0.5$, with $X$ and $Y$ as in Section~\ref{subsec:lls}, $W_0$ uniform on $\St(d,r)$, $M_0=I$ and $P_0=\Exp_I(0.8\,V)$ for a symmetric matrix $V$ with standard normal entries.
We take $R=120$, $K=150$, $\eta_x=0.2$, $\eta_y=0.2$, $s_{\max}=10$, $T=30$ for CG and NS, and $\lambda_0=10$, since at $\lambda_0=1$ the first outer step of RF$^2$SA overflows on some seeds.
Every lower-level update uses the exponential map of $\cS_{++}^d$, and the upper-level step uses the QR retraction on $\St(d,r)$.
In \RBDA{} the BB stepsize acts on the $M$ block, where the aggregated direction has both terms, whereas the $P$ block moves by the upper-level term $\mu\eta_y\ak\cG_P F$ alone.
The inner loop stops by the rule of Section~\ref{subsec:lls}, and the three quantities of Table~\ref{tab:nonlls} are recorded at the last outer step.

\begin{table}[ht]
\centering
\renewcommand{\arraystretch}{1.12}
\small
\setlength{\tabcolsep}{4pt}
\begin{tabular}{@{}lcccc@{}}
\toprule
\textbf{Method} & \textbf{Optimistic value $\vp(W)$} & $\rd(P,P^\star(W))$ & $\rd(M,M^\star(W))$ & \textbf{Time (s)} \\
\midrule
CG & $19.123\pm0.368$ & $5.80\pm0.71$ & $0.067\pm0.014$ & $7.1$ \\
NS & $19.109\pm0.371$ & $5.81\pm0.72$ & $0.068\pm0.014$ & $6.9$ \\
AD & $19.115\pm0.362$ & $5.79\pm0.71$ & $0.067\pm0.013$ & $14.7$ \\
RF$^2$SA~\cite{dutta2024riemannian} & $18.944\pm0.372$ & $\mathbf{(0.8\pm0.1)\times10^{-3}}$ & $(5.4\pm0.6)\times10^{-3}$ & $20.6$ \\
AdaRHD~\cite{shi2025adaptive} & $18.821\pm0.456$ & $5.33\pm0.63$ & $\mathbf{(7.0\pm6.4)\times10^{-5}}$ & $\mathbf{2.5}$ \\
\RBDA{} & $\mathbf{17.599\pm0.272}$ & $(1.5\pm1.0)\times10^{-3}$ & $0.025\pm0.004$ & $32.3$ \\
\bottomrule
\end{tabular}
\caption{Ignored-block problem, $20$ seeds, $R=120$: optimistic value, selection gap along the ignored block and transverse distance at the last outer step (mean $\pm$ standard deviation), and mean solver time.
Lower is better in every column, and bold marks the best entry per column.}
\label{tab:nonlls}
\end{table}

\begin{figure}[ht]
\centering
\includegraphics[width=\textwidth]{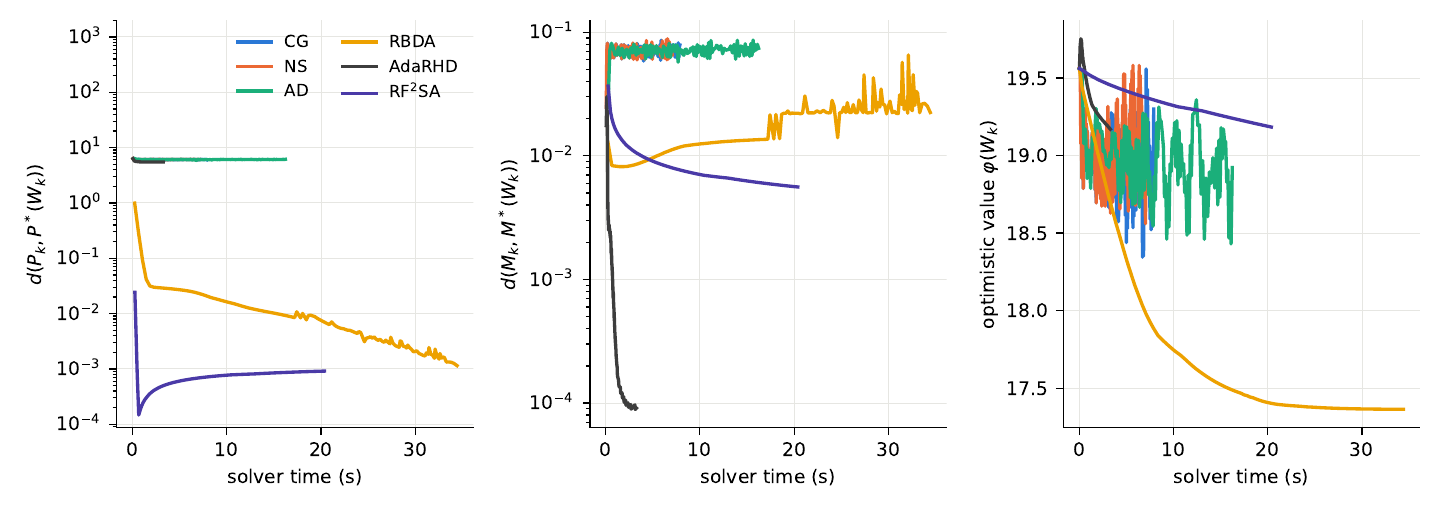}
\caption{Ignored-block problem, first seed, against solver time: (a) selection gap $\rd(P_k,P^\star(W_k))$, (b) transverse distance $\rd(M_k,M^\star(W_k))$, and (c) optimistic value $\vp(W_k)$.
The curves of CG, NS and AD coincide in panel (a).}
\label{fig:nonlls}
\end{figure}

Only the aggregated inner loop and RF$^2$SA move along the ignored block, as Table~\ref{tab:nonlls} and Figure~\ref{fig:nonlls} show.
\RBDA{} reduces the selection gap from $1.0$ after the first inner loop to $1.5\times10^{-3}$ at the last outer step and attains a smaller optimistic value than every baseline on all $20$ seeds, at four to five times the time of CG.
RF$^2$SA closes the selection gap as far as \RBDA{}, since its lower-level update also carries the upper-level gradient along $P$, and it still ends $1.3$ above \RBDA{} in $\vp$ on every seed, since its outer stepsize decreases with the multiplier and $R=120$ steps do not carry $W$ to the minimizer.
Selecting the right point of the solution set is necessary for the optimistic value and not sufficient.
In the transverse block \RBDA{} ends closer to $M^\star(W)$ than CG, NS and AD, since the upper-level term in~\eqref{A:agg} decays within each inner loop, whereas AdaRHD reaches the smallest transverse distance of the table and stays at the initial gap along $P$.
Appendix~\ref{app:tuning} reports the selection gap falling with $\mu$ on the swept grid.

\subsection{Few-shot meta-learning}\label{subsec:fewshot}
This experiment tests the selection on real data, where the episodes themselves underdetermine the lower level.
The upper-level variable is a shared representation $W\in\St(784,r)$, and the lower-level variable stacks one linear softmax head $\pi_i=(H_i,b_i)\in\br^{r\times5}\times\br^5$ per episode.
An episode $i$ draws five classes with one support and $15$ query images per class, its head is fit on the support features $W^\top u$, and the representation is scored by the query loss,
\begin{equation}\label{P:fs}
\min_{W\in\St(784,r)}\ \frac1E\sum_{i=1}^{E}\ell^{\mathrm{qry}}_i(W,\pi_i^\star)+\frac{\lambda}{2}\|\pi^\star\|^2
\quad\text{s.t.}\quad
\pi^\star\in\arg\min_{\pi}\ \frac1E\sum_{i=1}^{E}\ell^{\mathrm{sup}}_i(W,\pi_i),
\end{equation}
with $\ell^{\mathrm{sup}}_i$ and $\ell^{\mathrm{qry}}_i$ the cross-entropy of episode $i$ on its support and query features and $\lambda=10^{-4}$.
Each head fits five support points with $5r+5$ parameters, and five points in general position in $\br^r$ are linearly separable, so the support cross-entropy has infimum zero and no minimizer.
The lower-level solution set is empty, the problem lies outside the hypotheses of Section~\ref{sec:analysis}, and Appendix~\ref{app:fsreg} reports the variant with a ridge on the lower level.

We use $E=100$ training episodes as a full batch, $r=32$, $R=100$, $K=100$ run in full, $\eta_x=1$, $\eta_y=0.5$ and $\lambda_0=1$, with $W_0$ uniform on $\St(784,r)$ and a random initial head.
The lower level is Euclidean, so the inner update is additive and $s_{\max}=1/\hat L$, with $\hat L$ a power-iteration estimate of the largest eigenvalue of the lower-level Hessian at the initial head, and the upper-level step uses the QR retraction.
CG, NS and AdaRHD take $4K=400$ head steps before the implicit system, with $T=20$.
Since the lower-level Hessian $H$ is singular here, CG is also run on the Tikhonov-regularized system $(H+\varepsilon I)p=v$~\cite{hansen1998rank} with $\varepsilon=c\,\lambda_{\max}(H)$ and $c=10^{-2}$ fixed on the first seed, and this variant is labeled CG (Tikhonov).
A run is scored by meta-test accuracy on $50$ held-out episodes, on each of which the head is refit from a random start by $100$ gradient steps at the returned $W$ and evaluated on its query set.
The same protocol and the same $20$ episode draws are run on MNIST, Fashion-MNIST and KMNIST.
Table~\ref{tab:fs} also reports the terminal query loss on the training episodes, the upper-level objective of~\eqref{P:fs}.

\begin{table}[ht]
\centering
\renewcommand{\arraystretch}{1.12}
\small
\setlength{\tabcolsep}{4pt}
\begin{tabular}{@{}lccccc@{}}
\toprule
 & \multicolumn{3}{c}{\textbf{Meta-test accuracy}} & \shortstack{\textbf{Query loss}\\MNIST} & \\
\cmidrule(lr){2-4}
\textbf{Method} & \textbf{MNIST} & \textbf{Fashion-MNIST} & \textbf{KMNIST} & & \textbf{Time (s)} \\
\midrule
CG & $0.596\pm0.013$ & $0.660\pm0.022$ & $0.411\pm0.013$ & $3.467\pm0.196$ & $52.9$ \\
CG (Tikhonov) & $0.735\pm0.014$ & $0.659\pm0.016$ & $0.433\pm0.014$ & $0.437\pm0.026$ & $54.6$ \\
NS & $0.744\pm0.011$ & $0.668\pm0.018$ & $0.416\pm0.015$ & $0.459\pm0.015$ & $54.5$ \\
AD & $0.744\pm0.012$ & $0.679\pm0.016$ & $0.418\pm0.011$ & $0.677\pm0.043$ & $\mathbf{24.7}$ \\
RF$^2$SA~\cite{dutta2024riemannian} & $\mathbf{0.769\pm0.010}$ & $0.684\pm0.021$ & $\mathbf{0.483\pm0.014}$ & $0.305\pm0.010$ & $46.3$ \\
AdaRHD~\cite{shi2025adaptive} & $0.654\pm0.017$ & $0.609\pm0.016$ & $0.389\pm0.012$ & $1.092\pm0.037$ & $52.4$ \\
\RBDA{} & $0.747\pm0.012$ & $\mathbf{0.703\pm0.015}$ & $0.463\pm0.012$ & $\mathbf{0.212\pm0.005}$ & $79.2$ \\
\bottomrule
\end{tabular}
\caption{Few-shot meta-learning, $20$ episode draws, $R=100$: meta-test accuracy on three datasets (mean $\pm$ standard deviation, higher is better), terminal query loss on the training episodes on MNIST (lower is better), and mean solver time on MNIST.
Bold marks the best entry per column.}
\label{tab:fs}
\end{table}

\begin{figure}[ht]
\centering
\includegraphics[width=\textwidth]{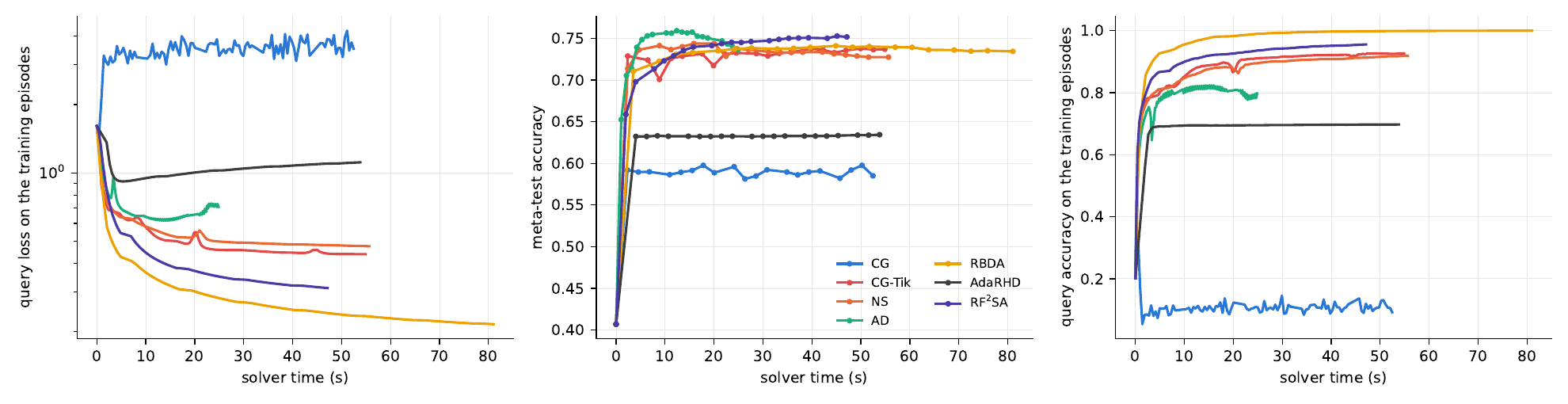}
\caption{Few-shot meta-learning on MNIST, first episode draw, against solver time: (a) query loss on the training episodes on a logarithmic axis, (b) meta-test accuracy, evaluated every five outer steps, and (c) query accuracy on the training episodes.}
\label{fig:fs}
\end{figure}

\RBDA{} attains the smallest query loss on all $20$ draws on each dataset, whereas in meta-test accuracy it is ahead of the four implicit and unrolled estimators on Fashion-MNIST and KMNIST and behind RF$^2$SA by two points on MNIST and KMNIST.
The query loss is a training quantity on $E=100$ episodes, and the head that minimizes it is not the one the held-out episodes reward most.
On MNIST the accuracy curves of NS, AD and CG (Tikhonov) peak within the first $40$ outer steps and decline by about two points afterwards.
Those of \RBDA{} and RF$^2$SA end within a point of their peaks.
Every method runs the inner cap in full here, so the extra time of \RBDA{} in Table~\ref{tab:fs} is a fixed cost per outer step, the second gradient in each inner iteration and the reverse sweep through the loop, and this section compares the returned heads.

CG fails because the linear system of~\eqref{A:HG1} is inconsistent here.
On five support points the lower-level Hessian of an episode has rank at most $20$, so $H$ is singular and the query gradient on the right-hand side lies largely in its null space.
Conjugate gradients diverge on such a system~\cite{kaasschieter1988preconditioned}, and on MNIST the query loss of CG rises from $1.61$ at $W_0$ to $3.47$.
NS is not affected, since its truncated Neumann series is a regularized inverse, and the Tikhonov shift brings CG (Tikhonov) to the level of NS at the price of the constant $c$ fixed on this problem.
AdaRHD solves the same singular system by conjugate gradients and inherits the divergence in a milder form.
RF$^2$SA stops at a query loss of $0.31$ against the $0.21$ of \RBDA{}, since its multiplier schedule sends the weight of the upper level in its lower-level update to zero.
On KMNIST the mean meta-test accuracy of every compared method is below one half.
Appendix~\ref{app:fsreg} adds a ridge to the lower level, under which CG recovers and the compared methods draw together, and Appendix~\ref{app:support} grows the support from one to ten samples per class.

\subsection{Data hyper-cleaning}\label{subsec:real}
This experiment is the data hyper-cleaning benchmark of~\cite{franceschi2017forward} in the configuration of~\cite[Sec.~8.2]{liu2022general}, the one problem on which BDA is also compared.
The split is $5{,}000$ training, $5{,}000$ validation and $60{,}000$ test images, with $50\%$ of the training labels replaced by a random wrong class.
The lower-level variable is a linear softmax head $\y=(W,b)\in\br^{784\times10}\times\br^{10}$, the upper-level variable a weight vector $\x\in\br^{5000}$ over the training samples, and
\begin{equation}\label{P:hc}
f(\x,\y)=\frac{\sum_{i\in\mathcal{D}_{\mathrm{tr}}}\sigma(\x_i)\,\ell(\y;u_i,v_i)}{\sum_{i\in\mathcal{D}_{\mathrm{tr}}}\sigma(\x_i)},
\qquad
F(\y)=\frac{1}{|\mathcal{D}_{\mathrm{val}}|}\sum_{i\in\mathcal{D}_{\mathrm{val}}}\ell(\y;u_i,v_i)+\frac{\lambda}{2}\|\y\|^2,
\end{equation}
with $\ell$ the cross-entropy, $\sigma$ the sigmoid and $\lambda=10^{-4}$.
The lower level carries no regularizer, so it is convex and nowhere strongly convex.
Its Hessian is singular along the shift $(W,b)\mapsto(W+v\mathbf{1}^\top,\,b+\beta\mathbf{1})$ and has near-zero eigenvalues along the directions of the head that the weighted training samples constrain only weakly, such as the weights of pixels that are almost always zero, so the LLS condition fails on real data.
The lower-level manifold is $\br^{7850}$, so Assumption~\ref{AN:geom} holds with $\kappa^-=0$, and $F$ is strongly convex in $\y$, so Assumption~\ref{AN:smooth} holds and Proposition~\ref{AN:sgrowth} gives the growth condition of Theorem~\ref{AN:rate} with $p=2$.
Since the data are not separable, Assumption~\ref{AN:qg} holds on bounded sets with a constant that we do not compute, and this problem is the one degenerate instance of this section within the hypotheses of Section~\ref{sec:analysis}.

We use $R=100$, $K=200$ run in full, $\eta_y=0.1$, $T=20$, $\lambda_0=1$, $s_{\max}=1/\hat L$ as in Section~\ref{subsec:fewshot}, and $\x_0=0$.
The upper-level step is the projected gradient step $\x\leftarrow P_{[-8,8]^{n}}\big(\x-\eta_x\widetilde{\cG}\vp(\x)\big)$ with $\eta_x=300$, which offsets the $1/n$ scale of the weighted loss.
CG, NS and AdaRHD take $4K=800$ head steps before the implicit system, and AdaRHD runs with $a_0=b_0=0.1$ on this problem in place of the $2$ used on the other three.
BDA~\cite{liu2022general} is also run, at the parameters of~\cite[Sec.~8.1]{liu2022general}, namely $\ak=1/(k+1)$, $\mu=0.5$ and the fixed stepsize $\eta_y$ in every inner iteration.
Beside the terminal test accuracy, Table~\ref{tab:hc-main} reports the peak of the test-accuracy curve over the run and the outer step at which it occurs, since several methods reach their best predictor within the first outer steps and lose it afterwards.

\begin{table}[ht]
\centering
\renewcommand{\arraystretch}{1.12}
\small
\setlength{\tabcolsep}{4pt}
\begin{tabular}{@{}lcccc@{}}
\toprule
\textbf{Method} & \textbf{Test accuracy} & \shortstack{\textbf{Peak test accuracy}\\(outer step)} & \textbf{Val.\ accuracy} & \textbf{Time (s)} \\
\midrule
CG & $0.8016\pm0.0045$ & $0.8208\pm0.0050$ ($1$) & $0.8538\pm0.0040$ & $219$ \\
NS & $0.7532\pm0.0044$ & $0.8208\pm0.0050$ ($1$) & $0.7599\pm0.0056$ & $223$ \\
AD & $0.8354\pm0.0035$ & $0.8378\pm0.0041$ ($49$) & $0.8445\pm0.0053$ & $\mathbf{102}$ \\
RF$^2$SA~\cite{dutta2024riemannian} & $0.8219\pm0.0044$ & $\mathbf{0.8774\pm0.0021}$ ($2$) & $\mathbf{0.8696\pm0.0038}$ & $161$ \\
AdaRHD~\cite{shi2025adaptive} & $0.5291\pm0.0057$ & $0.7617\pm0.0056$ ($1$) & $0.5338\pm0.0080$ & $218$ \\
BDA~\cite{liu2022general} & $0.8515\pm0.0029$ & $0.8516\pm0.0029$ ($98$) & $0.8610\pm0.0048$ & $216$ \\
\RBDA{} & $\mathbf{0.8565\pm0.0029}$ & $0.8565\pm0.0029$ ($99$) & $0.8694\pm0.0048$ & $219$ \\
\bottomrule
\end{tabular}
\caption{Hyper-cleaning on MNIST, $20$ seeds, $R=100$: terminal test accuracy, peak test accuracy along the run with the median outer step of the peak in parentheses, validation accuracy (mean $\pm$ standard deviation, higher is better), and mean solver time.
Bold marks the best entry per column.}
\label{tab:hc-main}
\end{table}

\begin{figure}[ht]
\centering
\includegraphics[width=\textwidth]{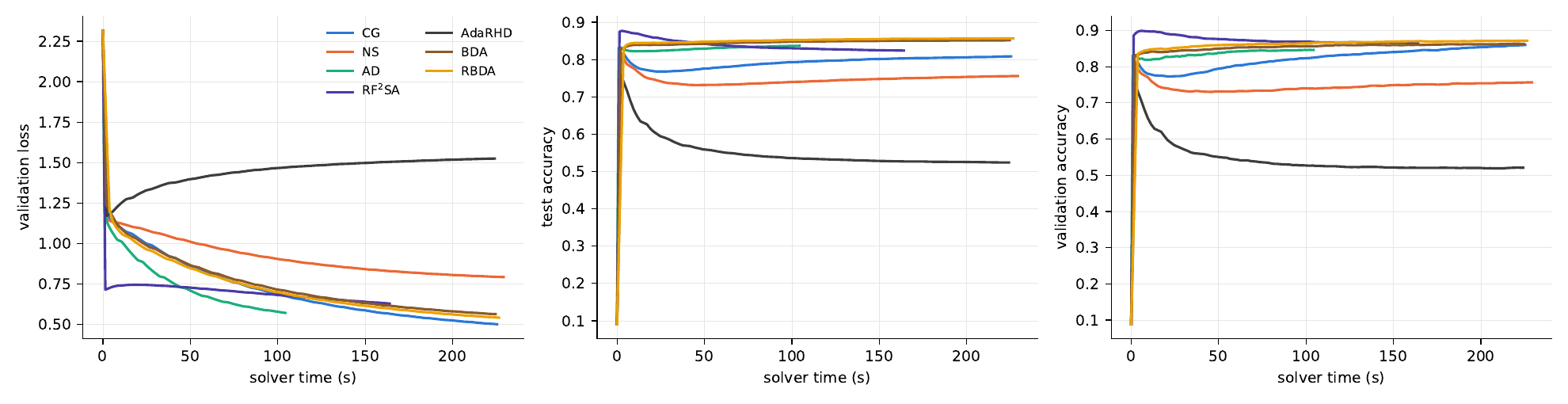}
\caption{Hyper-cleaning on MNIST, first seed, against solver time: (a) validation loss, (b) test accuracy, and (c) validation accuracy.}
\label{fig:hc-main}
\end{figure}

\RBDA{} attains the highest terminal test accuracy on all $20$ seeds against every baseline, by half a point over BDA and by two points over AD, and at the same time as BDA and the implicit estimators.
The implicit estimators lose accuracy during the first seconds of the run and do not recover it, so their peak in Table~\ref{tab:hc-main} is the accuracy after the first outer step.
The head they return is the $800$-step minimizer of the weighted training loss, and with no lower-level regularizer that head fits the corrupted half of the training set.
The validation accuracy of CG rises for the rest of the run and ends above that of AD, so its outer loop does optimize the validation objective.
What fails is the premise that an accurate lower-level solve and the implicit function theorem yield a predictor that generalizes.
NS adds the bias of its truncated series and is worse in both accuracies.
AdaRHD evaluates the same implicit formula at adaptive stepsizes, and its accuracy falls from that of its first head while its validation loss rises for the whole run.
RF$^2$SA reaches the highest peak of the table at the second outer step and then falls to a terminal value between AD and CG, since its multiplier schedule sends the weight of the upper level in its lower-level update to zero.
BDA is the closest competitor of \RBDA{}, and the two differ only in the decay exponent, the aggregation parameter and the BB stepsize, which Appendices~\ref{app:tuning} and~\ref{app:bb} vary.

Table~\ref{tab:hc} reports the terminal test accuracy of the same setting on Fashion-MNIST and KMNIST.
\RBDA{} is ahead of every baseline on all $20$ seeds on each dataset, by $1.5$ points over BDA on Fashion-MNIST and $1.1$ on KMNIST, and CG, NS and AdaRHD are the three lowest entries on each dataset.
On every dataset RF$^2$SA reaches the highest peak within the first outer steps and loses $5$ to $12$ points by the end, whereas \RBDA{} loses at most $2$.

\begin{table}[ht]
\centering
\renewcommand{\arraystretch}{1.12}
\small
\setlength{\tabcolsep}{5pt}
\begin{tabular}{@{}lcccc@{}}
\toprule
\textbf{Method} & \textbf{MNIST} & \textbf{Fashion-MNIST} & \textbf{KMNIST} & \textbf{Time (s)} \\
\midrule
CG & $0.8016\pm0.0045$ & $0.7204\pm0.0042$ & $0.5929\pm0.0049$ & $219$ \\
NS & $0.7532\pm0.0044$ & $0.6378\pm0.0057$ & $0.5445\pm0.0051$ & $223$ \\
AD & $0.8354\pm0.0035$ & $0.7616\pm0.0041$ & $0.6556\pm0.0041$ & $\mathbf{102}$ \\
RF$^2$SA~\cite{dutta2024riemannian} & $0.8219\pm0.0044$ & $0.7572\pm0.0029$ & $0.6414\pm0.0041$ & $161$ \\
AdaRHD~\cite{shi2025adaptive} & $0.5291\pm0.0057$ & $0.4826\pm0.0063$ & $0.3360\pm0.0046$ & $218$ \\
BDA~\cite{liu2022general} & $0.8515\pm0.0029$ & $0.7846\pm0.0033$ & $0.6850\pm0.0038$ & $216$ \\
\RBDA{} & $\mathbf{0.8565\pm0.0029}$ & $\mathbf{0.7997\pm0.0026}$ & $\mathbf{0.6956\pm0.0034}$ & $219$ \\
\bottomrule
\end{tabular}
\caption{Hyper-cleaning on three datasets, $20$ seeds, $R=100$: terminal test accuracy (mean $\pm$ standard deviation, higher is better) and mean solver time on MNIST.
Bold marks the best entry per column.}
\label{tab:hc}
\end{table}

\section{Conclusion}\label{sec:conclusion}
We proposed \RBDA{}, which extends bilevel descent aggregation to Riemannian manifolds for lower-level problems whose solution set is not a singleton.
Its inner loop aggregates the two descent directions in the tangent space under a decaying multiplier, and its hypergradient is the reverse-mode derivative of the unrolled loop.
For a Hadamard lower-level manifold and g-convex objectives with quadratic growth of the lower level, the inner iterates converge to a point of the optimistic solution set, at rates that depend on a growth condition of the upper-level objective.
Approximate minimizers of the objective with a finite number of inner iterations converge to minimizers of the optimistic value, and under the LLS condition the last inner iterate converges to the lower-level solution at the rate of the multiplier.

In the experiments the aggregation selects the optimistic solution where the implicit and unrolled estimators remain at the initial point of the ignored block.
It attains the highest test accuracy on the hyper-cleaning benchmark on three datasets, at two to three times the cost per outer step of the unrolled estimator.
The affine-invariant geometry keeps the lower-level problem tractable where the Euclidean descent of BDA reaches no solution.
The analysis covers the exponential map and the inner loop at a fixed upper-level variable, and a guarantee for the upper-level iteration of Algorithm~\ref{A:RBDA} and a last-iterate rate for the upper-level objective without the LLS condition remain open.
The BB stepsize is a heuristic that reduces the solver time on the synthetic problems, and the analysis does not cover the iterations that use it.

\section*{Declarations}

\paragraph{Funding}
The authors declare that no funds, grants, or other support were received during the preparation of this manuscript.

\paragraph{Competing interests}
The authors have no relevant financial or non-financial interests to disclose.

\paragraph{Data availability}
MNIST, Fashion-MNIST and KMNIST are loaded with \texttt{torchvision}, and the synthetic problem instances are generated by the released code from fixed random seeds.

\paragraph{Code availability}
\begin{sloppypar}
The code and the result files are available at \url{https://github.com/Cz1544252489/RBDA-Public}.
\end{sloppypar}

\bibliography{references}

\begin{appendices}
\section{Proofs for the Convergence Analysis}\label{app:proofs}

We collect here the proofs of the results in Section~\ref{sec:analysis}.
Table~\ref{tab:notation} lists the symbols of Section~\ref{sec:analysis} and of this appendix that are not in Table~\ref{tab:symbols}.

\begin{table}[H]
\centering\small
\renewcommand{\arraystretch}{1.2}
\begin{tabular}{@{}l >{\raggedright\arraybackslash}p{0.5\textwidth} >{\raggedright\arraybackslash}p{0.24\textwidth}@{}}
\toprule
\textbf{Symbol} & \textbf{Meaning} & \textbf{Introduced in} \\
\midrule
$\cS,\ \psi^*$ & solution set and minimum of $\psi$ at the fixed $\x$ & Section~\ref{subsec:setting} \\
$\kappa^-$ & lower bound on the sectional curvature of $\cN$ & Assumption~\ref{AN:geom} \\
$\by,\ \rho$ & optimistic solution used as center, and radius & \eqref{AN:ball} \\
$\bB,\ \bB_2$ & closed geodesic balls of radii $\rho$ and $2\rho$ about $\by$ & \eqref{AN:ball} \\
$\zeta$ & curvature factor of the comparison inequality & Lemma~\ref{AN:cmplemma} \\
$L_f,\ G_f,\ G_F$ & g-smoothness constant of $\psi$ on $\bB_2$, bounds on $\|\cG\psi\|$ and $\|\cG\Psi\|$ on $\bB$ & Assumption~\ref{AN:smooth} \\
$c$ & constant of the quadratic growth condition & Assumption~\ref{AN:qg} \\
$c',\ p$ & constant and order of the growth condition of $\Psi$ on $\cS$ & Theorem~\ref{AN:rate} \\
$C_f$ & strong g-convexity modulus of $\psi$ & Corollary~\ref{AN:lls} \\
$C_F$ & strong g-convexity modulus of $\Psi$ & Proposition~\ref{AN:sgrowth} \\
$\mu,\ \ak,\ s$ & aggregation parameter, multiplier, base stepsize & \eqref{A:agg}, \eqref{A:ss2} \\
$a,\ \theta$ & constants in the bounds $a(k+1)^{-\theta}\le\ak\le(k+1)^{-\theta}$ & Theorem~\ref{AN:rate} \\
$C_0$ & constant in the condition on $s$ & \eqref{AN:param} \\
$c_1,\ C_2$ & $s(1-\mu)$ and $2\zeta s^2\mu^2G_F^2$ & Lemma~\ref{AN:onestep} \\
$C_3$ & constant in the last-iterate bound & Lemma~\ref{AN:lastll} \\
$q$ & contraction factor under the LLS condition & Corollary~\ref{AN:lls} \\
$e_k,\ u_k$ & $\psi(\y_k)-\psi^*$ and $\Psi(\y_k)-\vp(\x)$ & \eqref{AN:ek} \\
$u_k^+,\ u_k^-$ & positive and negative parts of $u_k$ & \eqref{AN:uminus}, Appendix~\ref{app:pfmain} \\
$\Phi_k$ & squared distance from $\y_k$ to a comparison point & Lemma~\ref{AN:contain} \\
$j_k,\ B_0,\ C_4$ & index in $[k/2,k]$ and constants of the rates & Theorem~\ref{AN:rate}, Appendix~\ref{app:pfrate} \\
$\y_K(\x),\ \vp_K$ & $K$-th inner iterate at $\x$ and $F(\x,\y_K(\x))$ & Theorem~\ref{AN:consist} \\
$\cX,\ \varepsilon_K$ & compact set of upper-level variables and minimization tolerance & Theorem~\ref{AN:consist} \\
$\ty_j,\ \tilde\alpha_j,\ a_1$ & iterates and multipliers after the BB iterations, and the lower coefficient & Section~\ref{subsec:selection} \\
\bottomrule
\end{tabular}
\caption{Symbols of the analysis.}
\label{tab:notation}
\end{table}

\subsection{Auxiliary lemmas}\label{app:lemmas}

\begin{lemma}\label{AN:facts}
Let Assumptions~\ref{AN:geom} and~\ref{AN:smooth}\emph{(i)} hold, let $G_F:=\max_{\bB}\|\cG\Psi\|$, and let $\mu\in(0,1)$ and $\alpha\in(0,1]$.
Then the following hold:
\emph{(i)} if $\y\in\bB_2$, $v\in\cT_\y\cN$ and $\Exp_\y(v)\in\bB_2$, then $\psi(\Exp_\y(v))\le\psi(\y)+\langle\cG\psi(\y),v\rangle+\tfrac{L_f}{2}\|v\|^2$;
\emph{(ii)} $\|\cG\psi(\y)\|^2\le2L_f\big(\psi(\y)-\psi^*\big)$ for every $\y\in\bB$;
\emph{(iii)} $\|(1-\mu)\cG\psi(\y)+\mu\alpha\,\cG\Psi(\y)\|^2\le2(1-\mu)^2\|\cG\psi(\y)\|^2+2\mu^2\alpha^2G_F^2$ for every $\y\in\bB$.
\end{lemma}

\begin{proof}
For \emph{(i)}, the set $\bB_2$ is g-convex and $\psi$ is $L_f$-g-smooth on $\bB_2$.
On the Hadamard manifold $\cN$ we have $\Exp^{-1}_\y(\Exp_\y(v))=v$ and $\rd(\y,\Exp_\y(v))=\|v\|$, so \emph{(i)} is the quadratic upper bound in the definition of g-smoothness in Section~\ref{sec:preliminaries} at the pair $\y,\Exp_\y(v)$.

For \emph{(ii)}, let $\y\in\bB$ and $v:=-\cG\psi(\y)/L_f$.
Since $\|v\|\le G_f/L_f\le\rho$ by Assumption~\ref{AN:smooth}\emph{(i)}, we have $\rd(\Exp_\y(v),\by)\le\rd(\y,\by)+\|v\|\le2\rho$, and \emph{(i)} gives
$$
\psi^*\le\psi\big(\Exp_\y(v)\big)\le\psi(\y)+\langle\cG\psi(\y),v\rangle+\frac{L_f}{2}\|v\|^2=\psi(\y)-\frac{1}{2L_f}\|\cG\psi(\y)\|^2 .
$$

Statement \emph{(iii)} follows from $\|a+b\|^2\le2\|a\|^2+2\|b\|^2$ and $\|\cG\Psi(\y)\|\le G_F$ on $\bB$.
\end{proof}

For the iterates of~\eqref{AN:upd} we write
\begin{equation}\label{AN:ek}
e_k:=\psi(\y_k)-\psi^*\ge0,\qquad u_k:=\Psi(\y_k)-\vp(\x).
\end{equation}

\begin{lemma}\label{AN:onestep}
Let Assumptions~\ref{AN:geom} and~\ref{AN:smooth}\emph{(i)} hold, let $G_F:=\max_{\bB}\|\cG\Psi\|$, let $\mu\in(0,1)$, $\ak\in(0,1]$ and $0<s\le1/(4\zeta(1-\mu)L_f)$, and let $c_1:=s(1-\mu)$ and $C_2:=2\zeta s^2\mu^2G_F^2$.
If $\y_k\in\bB$ and $\y'\in\cS\cap\bB$, then
$$
\rd^2(\y_{k+1},\y')\le\rd^2(\y_k,\y')-c_1e_k+2s\mu\ak\big\langle\cG\Psi(\y_k),\Exp^{-1}_{\y_k}(\y')\big\rangle+C_2\ak^2 .
$$
\end{lemma}

\begin{proof}
By~\eqref{AN:upd} and the definition of $g_k$ in~\eqref{A:agg}, $\y_{k+1}=\Exp_{\y_k}(-s\,g_k)$.
Writing $\Phi_k:=\rd^2(\y_k,\y')$ and $E_k:=\Exp^{-1}_{\y_k}(\y')$, we have
\begin{align*}
\Phi_{k+1}
&\le\Phi_k+2s(1-\mu)\big\langle\cG\psi(\y_k),E_k\big\rangle+2s\mu\ak\big\langle\cG\Psi(\y_k),E_k\big\rangle+\zeta s^2\|g_k\|^2\\
&\le\Phi_k-2s(1-\mu)\,e_k+2s\mu\ak\big\langle\cG\Psi(\y_k),E_k\big\rangle+2\zeta s^2(1-\mu)^2\|\cG\psi(\y_k)\|^2+C_2\ak^2\\
&\le\Phi_k-2s(1-\mu)\big[1-2\zeta s(1-\mu)L_f\big]e_k+2s\mu\ak\big\langle\cG\Psi(\y_k),E_k\big\rangle+C_2\ak^2\\
&\le\Phi_k-c_1e_k+2s\mu\ak\big\langle\cG\Psi(\y_k),E_k\big\rangle+C_2\ak^2 .
\end{align*}
The first inequality is Lemma~\ref{AN:cmplemma} with $\y=\y_k$ and $v=-s\,g_k$.
The second follows from the first-order characterization of g-convexity of $\psi$ in Section~\ref{sec:preliminaries}, which gives $\langle\cG\psi(\y_k),E_k\rangle\le\psi(\y')-\psi(\y_k)=-e_k$, and from Lemma~\ref{AN:facts}\emph{(iii)}.
The third follows from Lemma~\ref{AN:facts}\emph{(ii)}, and the last from $s\le1/(4\zeta(1-\mu)L_f)$, which makes the bracket at least $\tfrac12$.
\end{proof}

\begin{lemma}\label{AN:contain}
Let Assumptions~\ref{AN:geom}--\ref{AN:radius} hold, let $\mu$, $\{\ak\}$ and $s$ satisfy~\eqref{AN:param}, and let $c_1$ and $C_2$ be as in Lemma~\ref{AN:onestep}, so that $C_0=C_2+2s^2\mu^2G_F^2c/c_1$.
Let $\y'\in\Sopt(\x)\cap\bB$ and $\Phi_k:=\rd^2(\y_k,\y')$.
Then for every $k\ge0$ the following hold:
\emph{(i)} $\Phi_{k+1}\le\Phi_k-c_1e_k-2s\mu\ak u_k+C_2\ak^2$;
\emph{(ii)} $\Phi_{k+1}\le\Phi_k-\tfrac{c_1}{2}e_k+C_0\ak^2$;
\emph{(iii)} $\y_k\in\bB$.
\end{lemma}

\begin{proof}
We first derive \emph{(i)} and \emph{(ii)} at an index $k$ with $\y_k\in\bB$, and then prove \emph{(iii)} by induction.

Let $\y_k\in\bB$ and $\y'\in\Sopt(\x)\cap\bB$.
The function $\Psi$ is g-convex on $\bB$ by Assumption~\ref{AN:smooth}\emph{(ii)}, so $\langle\cG\Psi(\y_k),\Exp^{-1}_{\y_k}(\y')\rangle\le\Psi(\y')-\Psi(\y_k)=-u_k$, and Lemma~\ref{AN:onestep} gives \emph{(i)}.

To derive \emph{(ii)}, let $\y_k^\cS:=P_\cS(\y_k)$.
Since $P_\cS(\by)=\by$ and $P_\cS$ is nonexpansive, $\rd(\y_k^\cS,\by)\le\rd(\y_k,\by)\le\rho$, so $\y_k^\cS\in\cS\cap\bB$.
Hence, with $u_k^-:=\max\{-u_k,0\}$,
\begin{equation}\label{AN:uminus}
-u_k\le\Psi(\y_k^\cS)-\Psi(\y_k)\le G_F\,\rd(\y_k,\cS)\le G_F\sqrt{c\,e_k},\qquad\text{so that}\qquad u_k^-\le G_F\sqrt{c\,e_k},
\end{equation}
where the first inequality follows from $\y_k^\cS\in\cS$ and the definition of $\vp(\x)$, the second from $\|\cG\Psi\|\le G_F$ along the geodesic from $\y_k$ to $\y_k^\cS$, which lies in the g-convex set $\bB$, and the third from Assumption~\ref{AN:qg}.
Combining \emph{(i)} with $-u_k\le u_k^-$ gives
\begin{align*}
\Phi_{k+1}
&\le\Phi_k-c_1e_k+2s\mu\ak u_k^-+C_2\ak^2\le\Phi_k-c_1e_k+2s\mu G_F\sqrt{c}\,\ak\sqrt{e_k}+C_2\ak^2\\
&\le\Phi_k-c_1e_k+\frac{c_1}{2}e_k+\frac{2s^2\mu^2G_F^2c}{c_1}\ak^2+C_2\ak^2=\Phi_k-\frac{c_1}{2}e_k+C_0\ak^2 ,
\end{align*}
where the second inequality follows from~\eqref{AN:uminus}, and the third from Young's inequality $2ab\le\epsilon a^2+b^2/\epsilon$ with $a=\sqrt{e_k}$, $b=s\mu G_F\sqrt{c}\,\ak$ and $\epsilon=c_1/2$.
This proves \emph{(ii)}.

We prove \emph{(iii)} by induction with $\y'=\by$.
The base case $\y_0\in\bB$ is Assumption~\ref{AN:radius}.
If $\y_0,\dots,\y_k\in\bB$, the inequality in \emph{(ii)} holds at every index $j\le k$, and summing it over $j=0,\dots,k$ gives
$$
\rd^2(\y_{k+1},\by)\le\rd^2(\y_0,\by)-\frac{c_1}{2}\sum_{j=0}^{k}e_j+C_0\sum_{j=0}^{k}\alpha_j^2\le\rd^2(\y_0,\by)+C_0\sum_{j\ge0}\alpha_j^2\le\rho^2 ,
$$
where the last inequality is~\eqref{AN:parama}.
Hence $\y_{k+1}\in\bB$, which completes the induction.
Since \emph{(i)}, \emph{(ii)} and~\eqref{AN:uminus} were derived under $\y_k\in\bB$ alone, they hold at every $k$.
\end{proof}

\subsection{Proof of Theorem~\ref{AN:main}}\label{app:pfmain}
We are now ready to prove Theorem~\ref{AN:main}.

\begin{proof}[Proof of Theorem~\ref{AN:main}]
By Lemma~\ref{AN:contain}\emph{(iii)}, $\y_k\in\bB$ for every $k$.
Summing Lemma~\ref{AN:contain}\emph{(ii)} with $\y'=\by$ over $k$ and using $\Phi_{k+1}\ge0$ and~\eqref{AN:parama} give
\begin{equation}\label{AN:esum}
\sum_{k\ge0}e_k\le\frac{2}{c_1}\Big(\rd^2(\y_0,\by)+C_0\sum_{k\ge0}\ak^2\Big)\le\frac{2\rho^2}{c_1},
\end{equation}
and in particular $e_k\to0$.

Let $u_k^+:=\max\{u_k,0\}$, so that $u_k^+=u_k+u_k^-$.
With $\y'=\by$ and $\Phi_k=\rd^2(\y_k,\by)$, we have
\begin{align*}
2s\mu\ak u_k^+
&=2s\mu\ak u_k+2s\mu\ak u_k^-\le\Phi_k-\Phi_{k+1}-c_1e_k+C_2\ak^2+2s\mu\ak u_k^-\\
&\le\Phi_k-\Phi_{k+1}-c_1e_k+C_2\ak^2+2s\mu G_F\sqrt{c}\,\ak\sqrt{e_k}\\
&\le\Phi_k-\Phi_{k+1}+\Big(C_2+\frac{s^2\mu^2G_F^2c}{c_1}\Big)\ak^2\le\Phi_k-\Phi_{k+1}+C_0\ak^2 ,
\end{align*}
where the first inequality is Lemma~\ref{AN:contain}\emph{(i)}, the second follows from~\eqref{AN:uminus}, the third from Young's inequality $2ab\le\epsilon a^2+b^2/\epsilon$ with $a=\sqrt{e_k}$, $b=s\mu G_F\sqrt{c}\,\ak$ and $\epsilon=c_1$, and the last from $C_0=C_2+2s^2\mu^2G_F^2c/c_1$.
Summing over $k$ and using $\Phi_{k+1}\ge0$ and~\eqref{AN:parama} give
\begin{equation}\label{AN:wsum}
\sum_{k\ge0}\ak u_k^+\le\frac{1}{2s\mu}\Big(\rd^2(\y_0,\by)+C_0\sum_{k\ge0}\ak^2\Big)\le\frac{\rho^2}{2s\mu}.
\end{equation}

We next show the following quasi-Fej\'er property of the iterates~\cite{combettes2001quasi}:
\begin{equation}\label{AN:fejer}
\lim_{k\to\infty}\rd^2(\y_k,\y')\ \text{exists for every }\y'\in\Sopt(\x)\cap\bB .
\end{equation}
Fix such a $\y'$ and let $\Phi_k:=\rd^2(\y_k,\y')$.
By Lemma~\ref{AN:contain}\emph{(ii)}, the sequence $R_k:=\Phi_k+C_0\sum_{j\ge k}\alpha_j^2$ satisfies
$$
0\le R_{k+1}=\Phi_{k+1}+C_0\sum_{j\ge k+1}\alpha_j^2\le\Phi_k-\frac{c_1}{2}e_k+C_0\ak^2+C_0\sum_{j\ge k+1}\alpha_j^2\le R_k ,
$$
so $\{R_k\}$ converges.
Since $\sum_{k\ge0}\ak^2<\infty$, the tail sum in $R_k$ tends to $0$, which proves~\eqref{AN:fejer}.

The bound~\eqref{AN:wsum} and $\sum_{k\ge0}\ak=\infty$ give $\liminf_{k\to\infty}u_k^+=0$.
Let $\{\y_{k_i}\}$ be a subsequence with $u_{k_i}^+\to0$.
Since $\bB$ is compact, we may assume, passing to a further subsequence, that $\y_{k_i}\to\hat\y$ for some $\hat\y\in\bB$.
Since $e_{k_i}\to0$, continuity of $\psi$ gives $\psi(\hat\y)=\psi^*$, that is, $\hat\y\in\cS$.
By~\eqref{AN:uminus}, $|u_{k_i}|=u_{k_i}^++u_{k_i}^-\le u_{k_i}^++G_F\sqrt{c\,e_{k_i}}\to0$, and continuity of $\Psi$ gives $\Psi(\hat\y)=\vp(\x)$.
Hence $\hat\y\in\Sopt(\x)\cap\bB$.

By~\eqref{AN:fejer} at $\y'=\hat\y$, the limit of $\rd^2(\y_k,\hat\y)$ exists, and it is $0$ along $\{k_i\}$, so $\y_k\to\hat\y$, which proves \emph{(ii)}.
Continuity of $\psi$ and $\Psi$ gives $\psi(\y_k)\to\psi(\hat\y)=\psi^*$ and $\Psi(\y_k)\to\Psi(\hat\y)=\vp(\x)$, which proves \emph{(i)}.
\end{proof}

\subsection{Proof of Theorem~\ref{AN:rate}}\label{app:pfrate}

\begin{proof}[Proof of Theorem~\ref{AN:rate}]
Since $\theta\in(\tfrac12,1)$, the bounds on $\ak$ give~\eqref{AN:paramb}, so the bounds~\eqref{AN:esum} and~\eqref{AN:wsum} hold.
Adding them and using $c_1=s(1-\mu)$ give
\begin{equation}\label{AN:B0}
\sum_{j\ge0}e_j+\sum_{j\ge0}\alpha_ju_j^+\le\frac{2\rho^2}{s(1-\mu)}+\frac{\rho^2}{2s\mu}=:B_0 .
\end{equation}

Fix $k\ge2$ and let $W_k:=\{j\in\mathbb{N}:k/2\le j\le k\}$, which contains at least $k/2$ indices.
For $j\in W_k$ we have $\alpha_j\ge a(j+1)^{-\theta}\ge a(k+1)^{-\theta}$, so~\eqref{AN:B0} gives
$$
\sum_{j\in W_k}\Big[e_j+a(k+1)^{-\theta}u_j^+\Big]\le B_0 .
$$
Hence some $j_k\in W_k$ satisfies
\begin{equation}\label{AN:goodindex}
e_{j_k}+a(k+1)^{-\theta}u_{j_k}^+\le\frac{2B_0}{k}.
\end{equation}
Both terms on the left of~\eqref{AN:goodindex} are nonnegative, so $e_{j_k}\le2B_0/k$, which is the first estimate in \emph{(i)}.
By Lemma~\ref{AN:contain}\emph{(iii)}, $\y_{j_k}\in\bB$, and Assumption~\ref{AN:qg} gives $\rd(\y_{j_k},\cS)\le\sqrt{c\,e_{j_k}}\le\sqrt{2cB_0/k}$, which is the second estimate in \emph{(i)}.

Since $(k+1)^\theta\le(2k)^\theta\le2k^\theta$ for $k\ge1$, the inequality~\eqref{AN:goodindex} also gives
\begin{equation}\label{AN:uplus}
u_{j_k}^+\le\frac{2B_0(k+1)^\theta}{a\,k}\le\frac{4B_0}{a\,k^{1-\theta}}.
\end{equation}
By~\eqref{AN:uminus} and \emph{(i)}, $u_{j_k}^-\le G_F\sqrt{c\,e_{j_k}}\le G_F\sqrt{2cB_0}\,k^{-1/2}\le G_F\sqrt{2cB_0}\,k^{-(1-\theta)}$, where the last inequality holds since $1-\theta<\tfrac12$ and $k\ge1$.
Hence $|u_{j_k}|=u_{j_k}^++u_{j_k}^-\le C_4k^{-(1-\theta)}$ with $C_4:=4B_0/a+G_F\sqrt{2cB_0}$, which is the first estimate in \emph{(ii)}.

Let $\y^\cS:=P_\cS(\y_{j_k})$.
Since $P_\cS(\by)=\by$ and $P_\cS$ is nonexpansive, $\rd(\y^\cS,\by)\le\rd(\y_{j_k},\by)\le\rho$, so $\y^\cS\in\cS\cap\bB$.
The geodesic from $\y_{j_k}$ to $\y^\cS$ lies in the g-convex set $\bB$, on which $\|\cG\Psi\|\le G_F$, so~\eqref{AN:uplus} and \emph{(i)} give
$$
\Psi(\y^\cS)-\vp(\x)\le u_{j_k}^++G_F\,\rd(\y_{j_k},\cS)\le\frac{4B_0}{a\,k^{1-\theta}}+\frac{G_F\sqrt{2cB_0}}{\sqrt{k}}\le\frac{C_4}{k^{1-\theta}},
$$
where the last inequality uses $k^{-1/2}\le k^{-(1-\theta)}$.
The growth condition on $\Psi$ at $\y^\cS\in\cS\cap\bB$ gives
$$
\rd\big(\y^\cS,\Sopt(\x)\big)\le\Big(\frac{C_4}{c'}\Big)^{1/p}k^{-(1-\theta)/p}.
$$
Since $\rd(\cdot,\Sopt(\x))$ is $1$-Lipschitz and $\rd(\y_{j_k},\y^\cS)=\rd(\y_{j_k},\cS)$, the triangle inequality and \emph{(i)} give
$$
\rd\big(\y_{j_k},\Sopt(\x)\big)\le\rd(\y_{j_k},\cS)+\rd\big(\y^\cS,\Sopt(\x)\big)\le\Big(\sqrt{2cB_0}+\Big(\frac{C_4}{c'}\Big)^{1/p}\Big)k^{-(1-\theta)/p},
$$
where the last inequality uses $k^{-1/2}\le k^{-(1-\theta)/p}$, which holds since $(1-\theta)/p<\tfrac12$.
This is the second estimate in \emph{(ii)}.
\end{proof}

\subsection{Proof of Proposition~\ref{AN:sgrowth}}\label{app:pfsgrowth}

\begin{proof}[Proof of Proposition~\ref{AN:sgrowth}]
Let $\y\in\cS\cap\bB$, and let $\gamma:[0,1]\to\cN$ be the geodesic with $\gamma(0)=\by$ and $\gamma(1)=\y$.
The sets $\cS$ and $\bB$ are g-convex, so $\gamma(t)\in\cS\cap\bB$ and $\Psi(\gamma(t))\ge\vp(\x)=\Psi(\by)$ for every $t\in[0,1]$.
Hence the right derivative of $\Psi\circ\gamma$ at $t=0$ is nonnegative, that is, $\langle\cG\Psi(\by),\Exp^{-1}_{\by}(\y)\rangle\ge0$.
The strong g-convexity inequality of Section~\ref{sec:preliminaries} at $\by$ gives
$$
\Psi(\y)-\vp(\x)\ge\big\langle\cG\Psi(\by),\Exp^{-1}_{\by}(\y)\big\rangle+\frac{C_F}{2}\rd^2(\by,\y)\ge\frac{C_F}{2}\rd^2(\by,\y)\ge\frac{C_F}{2}\rd^2\big(\y,\Sopt(\x)\big),
$$
where the last inequality follows from $\by\in\Sopt(\x)$.
If $\y\in\Sopt(\x)\cap\bB$, the left side is $0$, so $\y=\by$, which proves $\Sopt(\x)\cap\bB=\{\by\}$.
\end{proof}

\subsection{Proof of Theorem~\ref{AN:consist}}\label{app:pfconsist}
The proof relies on the following bound on the lower-level gap at the last iterate.

\begin{lemma}\label{AN:lastll}
Let Assumptions~\ref{AN:geom}--\ref{AN:radius} hold, let $\mu$, $\{\ak\}$ and $s$ satisfy~\eqref{AN:param}, and let $C_3:=s^2\mu^2G_F^2\big(L_f+2/(3c_1)\big)$.
Then for every $K\ge2$,
$$
\psi(\y_K)-\psi^*\le\frac{4\rho^2}{c_1K}+C_3\sum_{i\ge K/2}\alpha_i^2 .
$$
\end{lemma}

\begin{proof}
Let $k\ge0$ and $v_k:=-s\,g_k$, so that $\y_{k+1}=\Exp_{\y_k}(v_k)$.
By Lemma~\ref{AN:contain}\emph{(iii)}, $\y_k$ and $\y_{k+1}$ lie in $\bB\subset\bB_2$.
Writing $c_2:=\tfrac34c_1$, we have
\begin{align*}
e_{k+1}
&\le e_k+\big\langle\cG\psi(\y_k),v_k\big\rangle+\frac{L_f}{2}\|v_k\|^2\\
&\le e_k-s(1-\mu)\|\cG\psi(\y_k)\|^2+s\mu\ak G_F\|\cG\psi(\y_k)\|+L_fs^2(1-\mu)^2\|\cG\psi(\y_k)\|^2+L_fs^2\mu^2\ak^2G_F^2\\
&=e_k-s(1-\mu)\big[1-L_fs(1-\mu)\big]\|\cG\psi(\y_k)\|^2+s\mu\ak G_F\|\cG\psi(\y_k)\|+L_fs^2\mu^2G_F^2\ak^2\\
&\le e_k-c_2\|\cG\psi(\y_k)\|^2+s\mu\ak G_F\|\cG\psi(\y_k)\|+L_fs^2\mu^2G_F^2\ak^2\\
&\le e_k-\frac{c_2}{2}\|\cG\psi(\y_k)\|^2+C_3\ak^2\le e_k+C_3\ak^2 .
\end{align*}
The first inequality is Lemma~\ref{AN:facts}\emph{(i)} with $\y=\y_k$ and $v=v_k$.
The second follows from the Cauchy--Schwarz inequality, $\|\cG\Psi(\y_k)\|\le G_F$ and Lemma~\ref{AN:facts}\emph{(iii)}.
The third uses $L_fs(1-\mu)\le1/(4\zeta)\le\tfrac14$, which makes the bracket at least $\tfrac34$.
The fourth follows from Young's inequality $s\mu\ak G_F\|\cG\psi(\y_k)\|\le\tfrac{c_2}{2}\|\cG\psi(\y_k)\|^2+s^2\mu^2G_F^2\ak^2/(2c_2)$ and $1/(2c_2)=2/(3c_1)$.

Fix $K\ge2$ and let $W_K$ be the window of the proof of Theorem~\ref{AN:rate}, which contains at least $K/2$ indices.
Applying the last inequality from $j$ to $K-1$ gives $e_K\le e_j+C_3\sum_{i\ge K/2}\alpha_i^2$ for every $j\in W_K$.
Averaging over $j\in W_K$ and using~\eqref{AN:esum} give
$$
e_K\le\frac{1}{|W_K|}\sum_{j\in W_K}e_j+C_3\sum_{i\ge K/2}\alpha_i^2\le\frac{2}{K}\sum_{j\ge0}e_j+C_3\sum_{i\ge K/2}\alpha_i^2\le\frac{4\rho^2}{c_1K}+C_3\sum_{i\ge K/2}\alpha_i^2 .
$$
\end{proof}

\begin{proof}[Proof of Theorem~\ref{AN:consist}]
Let $f^*(\x):=\min_{\y\in\cN}f(\x,\y)$ and
$$
\tau_K:=\frac{4\rho^2}{c_1K}+C_3\sum_{i\ge K/2}\alpha_i^2 ,
$$
which does not depend on $\x\in\cX$ and tends to $0$ since $\sum_{i\ge0}\alpha_i^2<\infty$.
For every $\x\in\cX$ and $K\ge2$, Lemma~\ref{AN:lastll} at $\x$ gives $f(\x,\y_K(\x))-f^*(\x)\le\tau_K$.
By Lemma~\ref{AN:contain}\emph{(iii)} and Assumption~\ref{AN:radius} at $\x$, we have $\rd(\y_K(\x),\y_0)\le\rd(\y_K(\x),\by(\x))+\rd(\by(\x),\y_0)\le2\rho$, so every $\y_K(\x)$ lies in the compact set $\{\y\in\cN:\rd(\y,\y_0)\le2\rho\}$.
By Theorem~\ref{AN:main}\emph{(i)} at $\x$, we have $\vp_K(\x)=F(\x,\y_K(\x))\to\vp(\x)$ for every $\x\in\cX$.

Let $\bar\x$ be a limit point of $\{\x_K\}$, and let $\{\x_m\}_{m\in I}$ be a subsequence with $\x_m\to\bar\x$, where $I\subset\mathbb{N}$ is an infinite set of values of $K$, so that $\bar\x\in\cX$.
Passing to a further subsequence, we may assume that $\y_m(\x_m)\to\hat\y$ for some $\hat\y\in\cN$.
The function $f^*$ is upper semicontinuous as a pointwise infimum of the continuous functions $f(\cdot,\y)$, so continuity of $f$ gives
$$
0\le f(\bar\x,\hat\y)-f^*(\bar\x)\le\liminf_{m\to\infty}\Big[f\big(\x_m,\y_m(\x_m)\big)-f^*(\x_m)\Big]\le\lim_{m\to\infty}\tau_m=0,
$$
that is, $\hat\y\in\cS(\bar\x)$.
For every $\x\in\cX$, continuity of $F$ gives
\begin{equation}\label{AN:globchain}
\vp(\bar\x)\le F(\bar\x,\hat\y)=\lim_{m\to\infty}\vp_m(\x_m)\le\lim_{m\to\infty}\big(\vp_m(\x)+\varepsilon_m\big)=\vp(\x),
\end{equation}
where the first inequality follows from $\hat\y\in\cS(\bar\x)$ and the definition of $\vp$, and the second from the choice of $\x_m$.
This proves \emph{(i)}.

For \emph{(ii)}, every $\x\in\cX$ satisfies
$$
\limsup_{K\to\infty}\,\inf_{\cX}\vp_K\le\lim_{K\to\infty}\vp_K(\x)=\vp(\x),
$$
so $\limsup_{K\to\infty}\inf_{\cX}\vp_K\le\inf_{\cX}\vp$.
Choose a subsequence along which $\inf_{\cX}\vp_K$ converges to $\liminf_{K\to\infty}\inf_{\cX}\vp_K$ and, by compactness of $\cX$, along which $\x_K$ converges to some $\bar\x\in\cX$.
The argument for~\eqref{AN:globchain} applies along this subsequence, and $\vp_m(\x_m)\le\inf_{\cX}\vp_m+\varepsilon_m$ gives
$$
\inf_{\cX}\vp\le\vp(\bar\x)\le\lim_{m\to\infty}\vp_m(\x_m)\le\lim_{m\to\infty}\Big(\inf_{\cX}\vp_m+\varepsilon_m\Big)=\liminf_{K\to\infty}\inf_{\cX}\vp_K ,
$$
which proves \emph{(ii)}.
\end{proof}

\subsection{Proof of Corollary~\ref{AN:lls}}\label{app:pflls}

\begin{proof}[Proof of Corollary~\ref{AN:lls}]
Write $\y^*:=\y^*(\x)$, $\Phi_k:=\rd^2(\y_k,\y^*)$ and $e_k:=\psi(\y_k)-\psi^*$.
The strong g-convexity inequality of Section~\ref{sec:preliminaries} at $\y^*$, where $\cG\psi$ vanishes, gives
\begin{equation}\label{AN:scstar}
e_k\ge\frac{C_f}{2}\,\Phi_k .
\end{equation}
Assumption~\ref{AN:geom} holds with $\cS=\{\y^*\}$, so Lemmas~\ref{AN:facts} and~\ref{AN:onestep} apply.
For two distinct points $\y,\y'\in\bB_2$, the strong g-convexity inequality and Lemma~\ref{AN:facts}\emph{(i)} with $v=\Exp^{-1}_\y(\y')$ give $C_f\le L_f$.
Hence $s(1-\mu)C_f\le C_f/(4\zeta L_f)\le\tfrac14$, which gives $q\in[\tfrac{15}{16},1)$.

Let $\y_k\in\bB$, and let $c_1$ and $C_2$ be as in Lemma~\ref{AN:onestep}.
With $c=2/C_f$ we have $C_0=C_2+4s^2\mu^2G_F^2/(c_1C_f)$, and
\begin{align*}
\Phi_{k+1}
&\le\Phi_k-c_1e_k+2s\mu G_F\ak\sqrt{\Phi_k}+C_2\ak^2\le\Big(1-\frac{c_1C_f}{2}\Big)\Phi_k+2s\mu G_F\ak\sqrt{\Phi_k}+C_2\ak^2\\
&\le\Big(1-\frac{c_1C_f}{4}\Big)\Phi_k+\Big(C_2+\frac{4s^2\mu^2G_F^2}{c_1C_f}\Big)\ak^2=q\,\Phi_k+C_0\ak^2 .
\end{align*}
The first inequality is Lemma~\ref{AN:onestep} with $\y'=\y^*$, together with the Cauchy--Schwarz inequality, $\|\Exp^{-1}_{\y_k}(\y^*)\|=\sqrt{\Phi_k}$ and $\|\cG\Psi(\y_k)\|\le G_F$.
The second follows from~\eqref{AN:scstar}, and the third from Young's inequality $2ab\le\epsilon a^2+b^2/\epsilon$ with $a=\sqrt{\Phi_k}$, $b=s\mu G_F\ak$ and $\epsilon=c_1C_f/4$.
Since $q\le1$, the induction in the proof of Lemma~\ref{AN:contain}\emph{(iii)}, with this inequality in place of Lemma~\ref{AN:contain}\emph{(ii)}, shows that $\y_k\in\bB$ and $\Phi_{k+1}\le q\,\Phi_k+C_0\ak^2$ for every $k\ge0$.

Fix $k\ge2$.
Unrolling this inequality and splitting the sum at $k/2$ give
\begin{align*}
\Phi_k
&\le q^k\Phi_0+C_0\sum_{j=0}^{k-1}q^{\,k-1-j}\alpha_j^2
\le q^{k/2-1}\Big(\Phi_0+C_0\sum_{j<k/2}\alpha_j^2\Big)+C_0\max_{j\ge k/2}\alpha_j^2\sum_{j\ge k/2}q^{\,k-1-j}\\
&\le q^{k/2-1}\rho^2+\frac{C_0}{1-q}\max_{j\ge k/2}\alpha_j^2 ,
\end{align*}
where the second inequality uses $q^k\le q^{k/2-1}$ and $q^{\,k-1-j}\le q^{k/2-1}$ for $j<k/2$, and the last uses~\eqref{AN:parama} and $\sum_{j\le k-1}q^{\,k-1-j}\le1/(1-q)$.
This proves \emph{(i)}.

For \emph{(ii)}, every $j\ge k/2$ satisfies $j+1>k/2$, so the bound $\alpha_j\le a(j+1)^{-\theta}$ gives $\max_{j\ge k/2}\alpha_j^2\le a^2(k/2)^{-2\theta}=2^{2\theta}a^2k^{-2\theta}$.
Since $q<1$, the term $q^{k/2-1}\rho^2$ is $O(k^{-2\theta})$, and \emph{(i)} gives $\rd^2(\y_k,\y^*)=O(k^{-2\theta})$, which proves \emph{(ii)}.
\end{proof}

\section{Supplementary experiments}\label{app:supp}

\subsection{Sensitivity to the aggregation parameter}\label{app:tuning}
Table~\ref{tab:mu-lls} reports \RBDA{} on the singleton problem of Section~\ref{subsec:lls} on the first seed as the aggregation parameter $\mu$ varies, with every other setting as in that section.
At $\mu=0$ the inner loop is the plain lower-level descent with the BB stepsize.
On this grid the time to the $10^{-3}$ threshold, the inner iterations per outer step and the distance of the terminal inner iterate to $M^\star(W)$ are all smallest at $\mu=0$ and largest at $\mu=0.7$, whereas the outer steps to the threshold stay between $304$ and $307$.
On this problem the aggregation only adds to the cost of the inner loop, since the lower-level solution is unique and the upper-level gradient in~\eqref{A:agg} moves the inner iterates away from it.

\begin{table}[ht]
\centering
\renewcommand{\arraystretch}{1.12}
\small
\setlength{\tabcolsep}{3pt}
\begin{tabular}{@{}lcccccc@{}}
\toprule
$\mu$ & $0$ & $0.1$ & $0.2$ & $0.3$ & $0.5$ & $0.7$ \\
\midrule
time to $10^{-3}$ (s) & $36.1$ & $62.0$ & $66.3$ & $96.3$ & $104.1$ & $112.5$ \\
outer steps to $10^{-3}$ & $306$ & $306$ & $306$ & $307$ & $305$ & $304$ \\
inner iterations per outer step & $12.2$ & $23.9$ & $28.3$ & $28.4$ & $35.0$ & $50.2$ \\
$\rd(M_K,M^\star(W))$, last outer step & $5.2\times10^{-6}$ & $6.0\times10^{-3}$ & $8.0\times10^{-3}$ & $2.3\times10^{-2}$ & $7.1\times10^{-2}$ & $9.4\times10^{-2}$ \\
\bottomrule
\end{tabular}
\caption{Singleton problem, first seed, \RBDA{} with the aggregation parameter $\mu$ varied and every other setting as in Table~\ref{tab:lls}.
Lower is better in every row.}
\label{tab:mu-lls}
\end{table}

Table~\ref{tab:mu-nonlls} reports the same sweep on the ignored-block problem of Section~\ref{subsec:nonlls}.
The selection gap along the ignored block is smallest at $\mu=0.7$ and largest at $\mu=0.05$, since the $P$ block moves by the upper-level term alone, whereas the transverse distance is smallest at $\mu=0.05$ and stays below the $0.067$ of AD at every value of the grid.

\begin{table}[ht]
\centering
\renewcommand{\arraystretch}{1.12}
\small
\setlength{\tabcolsep}{3pt}
\begin{tabular}{@{}lcccccc@{}}
\toprule
$\mu$ & $0.05$ & $0.1$ & $0.2$ & $0.3$ & $0.5$ & $0.7$ \\
\midrule
$\rd(P,P^\star(W))$ & $7.5\times10^{-2}$ & $2.5\times10^{-2}$ & $8.8\times10^{-3}$ & $4.9\times10^{-3}$ & $2.8\times10^{-3}$ & $1.1\times10^{-3}$ \\
$\rd(M,M^\star(W))$ & $3.0\times10^{-4}$ & $6.3\times10^{-4}$ & $2.4\times10^{-3}$ & $4.3\times10^{-3}$ & $1.8\times10^{-2}$ & $2.2\times10^{-2}$ \\
$\vp(W)$ & $17.41$ & $17.38$ & $17.37$ & $17.37$ & $17.36$ & $17.36$ \\
\bottomrule
\end{tabular}
\caption{Ignored-block problem, first seed, \RBDA{} with the aggregation parameter $\mu$ varied and every other setting as in Table~\ref{tab:nonlls}, all three quantities at the last outer step.
Lower is better in every row.}
\label{tab:mu-nonlls}
\end{table}

Table~\ref{tab:mu-fs} reports the sweep on the few-shot problem of Section~\ref{subsec:fewshot} on MNIST.
The query loss on the training episodes is smallest at $\mu=0.9$, whereas the meta-test accuracy is largest at $\mu=0.7$, the fixed value, and lower at $\mu=0.9$.

\begin{table}[ht]
\centering
\renewcommand{\arraystretch}{1.12}
\small
\setlength{\tabcolsep}{3pt}
\begin{tabular}{@{}lccccccccc@{}}
\toprule
$\mu$ & $0.02$ & $0.05$ & $0.1$ & $0.15$ & $0.2$ & $0.3$ & $0.5$ & $0.7$ & $0.9$ \\
\midrule
meta-test accuracy & $0.659$ & $0.664$ & $0.714$ & $0.728$ & $0.724$ & $0.730$ & $0.729$ & $0.734$ & $0.725$ \\
query accuracy, training episodes & $0.514$ & $0.895$ & $0.987$ & $0.998$ & $0.999$ & $0.999$ & $1.000$ & $1.000$ & $1.000$ \\
query loss, training episodes & $4.60$ & $1.05$ & $0.57$ & $0.43$ & $0.36$ & $0.30$ & $0.24$ & $0.22$ & $0.19$ \\
\bottomrule
\end{tabular}
\caption{Few-shot meta-learning on MNIST, first episode draw, \RBDA{} with the aggregation parameter $\mu$ varied and every other setting as in Table~\ref{tab:fs}.
Higher is better in the first two rows and lower in the third.}
\label{tab:mu-fs}
\end{table}

Table~\ref{tab:mu-hc} reports the sweep on the hyper-cleaning problem of Section~\ref{subsec:real} on MNIST.
The validation loss, the upper-level objective, is smallest at $\mu=0.7$ and larger at $\mu=0.9$, whereas the test accuracy is largest at $\mu=0.9$.
At $\mu=0.1$ the test-accuracy curve peaks at the first outer step and falls afterwards, as the curves of the implicit estimators do in Section~\ref{subsec:real}, and from $\mu=0.3$ on it ends at its peak.
Table~\ref{tab:decay-hc} varies the decay exponent $\theta$ of~\eqref{A:ss2} on the same seed at the fixed $\mu$.
The slowest decay of the grid is best in both quantities, and $\theta=0.5$ lies outside the range $(\tfrac12,1]$ of Section~\ref{sec:analysis}.

\begin{table}[ht]
\centering
\renewcommand{\arraystretch}{1.12}
\small
\setlength{\tabcolsep}{5pt}
\begin{tabular}{@{}lccccc@{}}
\toprule
$\mu$ & $0.1$ & $0.3$ & $0.5$ & $0.7$ & $0.9$ \\
\midrule
test accuracy & $0.829$ & $0.836$ & $0.844$ & $0.857$ & $0.875$ \\
peak test accuracy (outer step) & $0.832$ ($1$) & $0.836$ ($98$) & $0.844$ ($100$) & $0.857$ ($99$) & $0.875$ ($99$) \\
validation loss & $0.575$ & $0.562$ & $0.549$ & $0.542$ & $0.562$ \\
\bottomrule
\end{tabular}
\caption{Hyper-cleaning on MNIST, first seed, \RBDA{} with the aggregation parameter $\mu$ varied and every other setting as in Table~\ref{tab:hc-main}.
Higher is better in the first two rows and lower in the third.}
\label{tab:mu-hc}
\end{table}

\begin{table}[ht]
\centering
\renewcommand{\arraystretch}{1.12}
\small
\setlength{\tabcolsep}{6pt}
\begin{tabular}{@{}lccc@{}}
\toprule
$\theta$ & $0.5$ & $0.75$ & $1$ \\
\midrule
test accuracy & $0.863$ & $0.857$ & $0.853$ \\
validation loss & $0.509$ & $0.542$ & $0.559$ \\
\bottomrule
\end{tabular}
\caption{Hyper-cleaning on MNIST, first seed, \RBDA{} with the decay exponent $\theta$ varied and every other setting as in Table~\ref{tab:hc-main}.
Higher is better in the first row and lower in the second.}
\label{tab:decay-hc}
\end{table}

\subsection{Regularizing the few-shot lower level}\label{app:fsreg}
Table~\ref{tab:fs-reg} adds a ridge $\tfrac{\lambda_\ell}{2}\|\pi\|^2$ to the lower level of~\eqref{P:fs} on the first episode draw, so that for $\lambda_\ell>0$ the lower-level solution is unique and the problem is in the LLS setting.
The ridge makes the implicit system nonsingular, and CG, which is fifteen points below the other methods at $\lambda_\ell=0$, is within a point of CG (Tikhonov) from $\lambda_\ell=10^{-3}$ on.
At $\lambda_\ell=10^{-1}$ every method is lower than at $10^{-2}$, since the ridge then dominates the support loss.
The table reports one draw, whose accuracy fluctuates by about $0.014$ between draws in Section~\ref{subsec:fewshot}, so differences of one or two points between the methods at a fixed $\lambda_\ell$ are within that fluctuation.

\begin{table}[ht]
\centering
\renewcommand{\arraystretch}{1.12}
\small
\setlength{\tabcolsep}{5pt}
\begin{tabular}{@{}lccccc@{}}
\toprule
$\lambda_\ell$ & $0$ & $10^{-4}$ & $10^{-3}$ & $10^{-2}$ & $10^{-1}$ \\
\midrule
CG & $0.585$ & $0.668$ & $0.771$ & $0.749$ & $\mathbf{0.697}$ \\
CG (Tikhonov) & $0.737$ & $\mathbf{0.767}$ & $0.773$ & $0.748$ & $\mathbf{0.697}$ \\
NS & $0.728$ & $\mathbf{0.767}$ & $0.773$ & $0.749$ & $\mathbf{0.697}$ \\
AD & $0.736$ & $0.765$ & $0.770$ & $0.743$ & $0.697$ \\
RF$^2$SA & $\mathbf{0.752}$ & $0.757$ & $0.765$ & $0.734$ & $0.653$ \\
AdaRHD & $0.634$ & $0.748$ & $0.765$ & $0.738$ & $0.672$ \\
\RBDA{} & $0.734$ & $0.745$ & $\mathbf{0.775}$ & $\mathbf{0.750}$ & $\mathbf{0.697}$ \\
\bottomrule
\end{tabular}
\caption{Few-shot meta-learning on MNIST, first episode draw, with a ridge $\lambda_\ell$ on the lower level and every other setting as in Table~\ref{tab:fs}: meta-test accuracy, higher is better.
Bold marks the best entry per column, with ties at three decimals resolved by the unrounded values.}
\label{tab:fs-reg}
\end{table}

The constant $c=10^{-2}$ of the Tikhonov shift in Section~\ref{subsec:fewshot} is the value among $10^{-4}$, $10^{-3}$, $10^{-2}$ and $10^{-1}$ at which CG (Tikhonov) attains its smallest query loss on the first draw.

\subsection{Support size in few-shot meta-learning}\label{app:support}
Table~\ref{tab:kshot} grows the support of the few-shot problem from one to ten samples per class on the first draw, with the query set and every other setting as in Section~\ref{subsec:fewshot}.
With more support samples each head is determined more closely by its support set, and the room for the selection shrinks.
At $k=5$ and $k=10$ the five methods lie within $1.3$ points of one another in meta-test accuracy, against $2.4$ points at $k=1$.
\RBDA{} has the smallest query loss at $k\in\{1,2,5\}$ and RF$^2$SA at $k=10$, whereas the accuracy is highest for RF$^2$SA at $k=1$ and for AD at $k\in\{2,5,10\}$.

\begin{table}[ht]
\centering
\renewcommand{\arraystretch}{1.12}
\small
\setlength{\tabcolsep}{5pt}
\begin{tabular}{@{}lcccccccc@{}}
\toprule
 & \multicolumn{4}{c}{\textbf{Meta-test accuracy}} & \multicolumn{4}{c}{\textbf{Query loss}} \\
\cmidrule(lr){2-5}\cmidrule(lr){6-9}
\textbf{Method} & $k=1$ & $2$ & $5$ & $10$ & $k=1$ & $2$ & $5$ & $10$ \\
\midrule
CG (Tikhonov) & $0.737$ & $0.801$ & $0.868$ & $0.900$ & $0.44$ & $0.46$ & $0.42$ & $0.44$ \\
NS & $0.728$ & $0.818$ & $0.873$ & $0.896$ & $0.48$ & $0.43$ & $0.43$ & $0.45$ \\
AD & $0.736$ & $\mathbf{0.829}$ & $\mathbf{0.874}$ & $\mathbf{0.909}$ & $0.71$ & $0.51$ & $0.39$ & $0.36$ \\
RF$^2$SA & $\mathbf{0.752}$ & $\mathbf{0.829}$ & $0.868$ & $0.903$ & $0.31$ & $0.32$ & $0.31$ & $\mathbf{0.31}$ \\
\RBDA{} & $0.734$ & $0.813$ & $0.866$ & $0.898$ & $\mathbf{0.22}$ & $\mathbf{0.24}$ & $\mathbf{0.29}$ & $0.35$ \\
\bottomrule
\end{tabular}
\caption{Few-shot meta-learning on MNIST, first episode draw, with $k$ support samples per class and every other setting as in Table~\ref{tab:fs}.
Higher is better for the accuracy and lower for the loss, and bold marks the best entry per column, with ties at three decimals resolved by the unrounded values.}
\label{tab:kshot}
\end{table}

\subsection{The BB stepsize}\label{app:bb}
Table~\ref{tab:bb} compares \RBDA{} with \RBDAnoBB{}, which is Algorithm~\ref{A:RBDA} with $K_1=0$, on the four problems of Section~\ref{sec:experiments} with every other setting unchanged.
On the two synthetic problems, where the inner loop stops early, the BB stepsize reduces the inner iterations per outer step and the solver time, and on the singleton problem it leaves the outer trajectory unchanged.
On the ignored-block problem it enlarges the selection gap and the transverse distance at an unchanged optimistic value.
On the two real-data problems, where the inner loop runs in full, it changes the returned head at the same time.
On few-shot meta-learning it halves the query loss on MNIST, and \RBDAnoBB{} has the higher meta-test accuracy on MNIST and Fashion-MNIST.
On hyper-cleaning \RBDAnoBB{} has the higher test accuracy on MNIST and KMNIST and \RBDA{} on Fashion-MNIST, within a point in each case.
On the singleton problem the fraction $K_1/K$ of BB iterations was also varied over $\{0,\tfrac14,\tfrac12,\tfrac34,1\}$ on the first seed.
The time to the $10^{-3}$ threshold is $89.0$~s at $0$ and $52.6$~s at $1$, and the distance of the terminal inner iterate to $M^\star(W)$ is $6.1\times10^{-2}$ up to the fraction $\tfrac12$ and about $10^{-1}$ at $\tfrac34$ and $1$.

\begin{table}[ht]
\centering
\renewcommand{\arraystretch}{1.12}
\small
\setlength{\tabcolsep}{5pt}
\begin{tabular}{@{}llcc@{}}
\toprule
\textbf{Problem} & \textbf{Quantity} & \RBDA{} & \RBDAnoBB{} \\
\midrule
Singleton & time to $10^{-3}$ (s) & $73.2\pm7.6$ & $92.2\pm7.0$ \\
 & inner iterations per outer step & $47.7$ & $80.0$ \\
Ignored block & $\vp(W)$ & $17.599\pm0.272$ & $17.599\pm0.272$ \\
 & $\rd(P,P^\star(W))$ & $(1.5\pm1.0)\times10^{-3}$ & $(0.8\pm0.6)\times10^{-3}$ \\
 & $\rd(M,M^\star(W))$ & $0.025\pm0.004$ & $0.012\pm0.001$ \\
 & time (s) & $32.3$ & $47.7$ \\
Few-shot & meta-test accuracy, MNIST & $0.747\pm0.012$ & $0.766\pm0.010$ \\
 & meta-test accuracy, Fashion-MNIST & $0.703\pm0.015$ & $0.711\pm0.016$ \\
 & meta-test accuracy, KMNIST & $0.463\pm0.012$ & $0.462\pm0.012$ \\
 & query loss, MNIST & $0.212\pm0.005$ & $0.426\pm0.014$ \\
 & time (s) & $79.2$ & $78.4$ \\
Hyper-cleaning & test accuracy, MNIST & $0.8565\pm0.0029$ & $0.8601\pm0.0028$ \\
 & test accuracy, Fashion-MNIST & $0.7997\pm0.0026$ & $0.7978\pm0.0026$ \\
 & test accuracy, KMNIST & $0.6956\pm0.0034$ & $0.7024\pm0.0033$ \\
 & time (s) & $219$ & $216$ \\
\bottomrule
\end{tabular}
\caption{\RBDA{} and \RBDAnoBB{} on the four problems, $20$ seeds each, with the settings of Tables~\ref{tab:lls}, \ref{tab:nonlls}, \ref{tab:fs} and~\ref{tab:hc}.
Higher is better for the accuracies and lower for every other quantity.}
\label{tab:bb}
\end{table}

\end{appendices}

\end{document}